\documentclass[10pt]{amsart}
\usepackage[cp1251]{inputenc}
\usepackage[english,russian]{babel}
\usepackage{amsmath}
\usepackage{amssymb}
\usepackage{amsfonts}
\usepackage{bm}%
\usepackage{amscd}%
\usepackage{color}%
\usepackage{xcolor}
\usepackage{xcolor}
\usepackage{hyperref}

\newtheorem{lemma}{Lemma}
\newtheorem{theorem}{Theorem}
\newtheorem{definition}{Definition}

\newtheorem{proposition}{Proposition}

\begin{document}
\renewcommand{\refname}{References}

\thispagestyle{empty}

\title[Funk--Minkowski  transform and  spherical convolution of Hilbert type
   % in reconstruction of functions on the sphere
] 
       {Funk--Minkowski  transform and  spherical convolution of Hilbert type in reconstructing functions on the sphere}
\author{{S.\,G. Kazantsev}}%
\address{Sergei Gavrilovich  Kazantsev 
\newline\hphantom{iii} Sobolev Institute of Mathematics,
\newline\hphantom{iii} pr. Koptyuga, 4,
\newline\hphantom{iii} 630090, Novosibirsk, Russia}%
\email{kazan@math.nsc.ru}%

\thanks{\sc Kazantsev, S.\,G., Funk--Minkowski  transform and  spherical convolution of Hilbert type    in reconstructing functions on the sphere}
\thanks{\copyright \ 2018 Kazantsev S.G.}

 \vspace{1cm} 
\maketitle
 {\small
\begin{quote}
\noindent{\sc Abstract. } The Funk--Minkowski transform
 ${\mathcal F}$  associates a function $f$ on the sphere ${\mathbb S}^2$
  with its mean  values (integrals) along all great circles of the sphere.  Thepresented  analytical inversion formula  reconstruct the unknown function $f$ completely if  two Funk--Minkowski transforms, ${\mathcal F}f$ and ${\mathcal F} \nabla f$, are known. Another result  of this
 article is related to the problem of  Helmholtz--Hodge decomposition for
 tangent vector field on the sphere ${\mathbb S}^2$. We proposed solution for this problem which is used  the Funk--Minkowski  transform ${\mathcal F}$ and  Hilbert type spherical convolution ${\mathcal S}$.
\medskip

 \end{quote}
}

%MSC: primary 44A12; secondary 52A38
%Math Subject Classifications. 65T40, 45Q05, 65N21, 44A12, 62G05.

\section{Introduction}

 The paper is devoted to the  analytical  inverse of the
 Minkowski--Funk  transform (F--M transform). This transform  was introduced by
 P. Funk ~\cite{Funk:1913, Funk:1915, Funk:1916},
 based on the work ~\cite{Minkowski:1904} of  H. Minkowski.
 In literature  Funk--Minkowski transform is  known also
 as the Funk transform, Funk--Radon transform or spherical Radon
 transform. F--M transform associates a function on
 the unit sphere ${\mathbb S}^2$ in
 ${\mathbb R}^3$ with its mean  values (integrals) along all great circles of the sphere.
  Funk--Minkowski transform is a geodesic transform because the
  great circles on the sphere are geodesics. In recent time many
  authors investigate the generalized  Funk--Minkowski transforms
  (or nongeodesic Funk--Minkowski transforms) on the sphere ${\mathbb S}^2$, which
  include nongeodesic paths of integration, such as
  circles with fixed diameter ~\cite{Schneider:1969, Rubin:2000, NattererWubelling:2001},
  circles perpendicular to the equator
  ~\cite{ZangerlScherzer:2010, HielscherQuellmalz:2016, GindikinReedsShepp:1994, NattererWubelling:2001}
  and  circles, which   obtained
  by intersections of  the sphere with planes passing through
  a fixed common point ${\bf a} \in {\mathbb R}^3$, for example,  through the  northpole
  ${\bf k} \in {\mathbb S}^2$
  ~\cite{Salman:2016, Quellmalz:2017, AbouelazDaher:1993,
   Daher:2001, Helgason:2011, Palamodov:2017}.

 Funk--Minkowski transform
 plays an important role in the study of other
 integral
 transforms on the sphere and has various applications, for example,
 it is used in the convex geometry, harmonic analysis,
 image processing  and in photoacoustic tomography, see ~\cite{Minkowski:1904,
  Stepanov:2017, Dann:2010,
  LouisRiplingerSpiessSpodarev:2011, Louis:2016,
  QuellmalzHielscherLouis:2018, YarmanYazici:2011,
  ZangerlScherzer:2010, HristovaaMoonSteinhauer:2016}.

  Let ${\mathbb B}^3$ and  ${\mathbb S}^2$ be the unit ball and
  the unit sphere in $\mathbb R^3$, respectively,
  i.e.
  ${\mathbb B}^3 = \{{\bf x} \in \mathbb R^3 : |{\bf x}|<1 \}$ and
  ${\mathbb S}^2 = \partial \mathbb B^3 =\{{\bm \xi} \in \mathbb R^3 : |{\bm \xi}|=1 \},$
  where $| \cdot|$ denotes the Euclidean norm.
  Throughout the paper we adopt the convention to denote in bold type the vectors in
  ${\mathbb R}^3,$ and in simple type the scalars in ${\mathbb R}.$ By the  greek
   letters ${\bm \theta}$, ${\bm \eta}$, ${\bm \xi}$ and so on we denote the units
   vectors ${\mathbb S}^2$. We will use for unit vector ${\bm \xi}$ on the sphere $\mathbb{S}^2$
  usual angular coordinates $(\theta, \varphi)$
  \[
  {\bm\xi}={\bm\xi}({\theta},{\varphi})={\bf i}\sin{\theta}\cos{\varphi}+
  {\bf  j}\sin{\theta}\sin{\varphi}+ {\bf  k}\cos{\theta}=
  ({\sin\theta} \cos {\varphi},  {\sin\theta} \sin {\varphi},\cos {\theta}),
  \]
  where   $0 < \theta < \pi$ (the colatitude),
  $0 < \varphi < 2 \pi$  (the longitude) and
  %$\theta$ --- latitude,
  $t=\cos \theta$ --- polar
  distance.

  The  plane ${\bm\xi}^\perp=\{{\bf x}\in\mathbb{R}^3 : {\bf x}\centerdot {\bm\xi}=0\}$
  is spanned by the two orthonormal
  vectors ${\bf e}_{1},$ ${\bf e}_{2}$  with
  representations in polar coordinates
  \[
   {\bf e}_{1}(\bm\xi)=\frac{\partial{\bm\xi}}{\partial {\theta}}=
  ({\cos\theta}  \cos {\varphi}, {\cos\theta}  \sin {\varphi},-\sin {\theta}),
  \,
  {\bf e}_2(\bm\xi)=\frac{1}{\sin{\theta}}\frac{\partial{\bm\xi}}{\partial {\varphi}}=
  ( -\sin {\varphi}, \cos {\varphi}, 0).
   \]
   The vectors ${\bf e}_{1}({\bm\xi}),$
  ${\bf e}_{2}({\bm\xi})$ and   ${\bm\xi}$
  form the so called local moving
  triad ${\bm\xi} \centerdot  {\bf e}_{1}=0$,
  ${\bm\xi} \centerdot  {\bf e}_{2}=0$,
  ${\bf e}_{1}\centerdot  {\bf e}_{2}=0$,
  where  $ \centerdot $ denotes the inner product
  of two vectors in ${\mathbb R}^3$.

 Let  denote by $f_{even}$ and
 $f_{odd}$ the even and odd parts of function $f$ on ${\mathbb S}^2$, respectively, that is,
 we have
 \[
 f({\bm \xi})=f_{even}({\bm \xi})+f_{odd}({\bm
  \xi}), \
 f_{even}({\bm \xi})=\frac{f({\bm \xi})+f(-{\bm \xi})}{2}, \ \
 f_{odd}({\bm \xi})=
 \frac{f({\bm \xi})-f(-{\bm \xi})}{2}.
 \]
%In the following, N^{+}=\{ 1,2, ..}}$ is the set of all
% $  non-zeronatural numbers, $N^{+}=N+\{0 \}$.

 The space of continuous functions on the sphere ${\mathbb S}^2$ is
 denoted by $C({\mathbb S}^2)$  and is  endowed with the supremum
 norms
 \[
 ||f||_{C({\mathbb S}^2)}=
 {\rm sup}_{{\bm \xi}\in{\mathbb S}^2 }|f({\bm \xi})|.
 \]
 $C({\mathbb S}^2)$, $C_{even}({\mathbb S}^2)$ and $C_{odd}({\mathbb S}^2)$
 denote the space of continuous functions on ${\mathbb S}^2$,
  the space of even  continuous functions on ${\mathbb S}^2$
  and the space of odd  continuous functions on ${\mathbb S}^2$,
 respectively. The subset of  $C_{even}({\mathbb S}^2)$ ($C_{odd}({\mathbb S}^2)$) that
 contains the infinitely differentiable functions will be denoted
 by $C^{\infty}_{even}({\mathbb S}^2)$ ($C^{\infty}_{odd}({\mathbb S}^2)$).

 \begin{definition} Let $f$ be a
 continuous function on the sphere ${\mathbb S}^2$,
 $f \in  C({\mathbb S}^2)$.
 Then, for a unit vector  ${\bm \xi} \in {\mathbb S}^2$
 the Funk--Minkowski transform  of a function  $f$ is a
  function ${\mathcal F}f$ on ${\mathbb S}^2$,
  given by
    \begin{align}
   \label{SEMR:Funk--Minkowski-scalar}
  \{ {\mathcal F}f \}({\bm \xi})\equiv{\mathcal F}_{\bm \xi}f=
\frac{1}{2\pi}
  \int_{0}^{2\pi}f \Big ({\bf e}_{1}({\bm \xi})\cos \omega
 +{\bf e}_{2}({\bm \xi})\sin \omega \Big ) \,{\rm d}{\omega} .
  \end{align}
 \end{definition}
It is clear that the Funk--Minkowski transform is even,
 $\{ {\mathcal F}f\}(-{\bm \xi})=\{ {\mathcal F}f\}({\bm \xi})$, and
${\mathcal F}$ annihilates all odd functions.

  The inversion of the Funk--Minkowski transform has been treated by
  many authors and there are exist several inversion formulas
  in the literature, see
  ~\cite{Funk:1913, Semyanistyi:1961, Helgason:1999, Rubin:2002,  Rubin:2003}.
  In ~\cite{Funk:1913, Funk:1916} P. Funk
     proved that an even  function  can be recovered from the knowledge
  of integrals over great circles and
    presented  two different inversion methods: the first method is
  based on the spherical harmonic decomposition of the functions
  $f$,  ${\mathcal F}f$ and the second one utilizes Abel's integral equation,
  ~\cite{NattererWubelling:2001}.
 %and the only even part of a function $f$ over ${\mathbb S}^2$
 %can be uniquely determined from its transform
 % (\ref{SEMR:Funk--Minkowski-scalar}).

The inversion formula after P. Funk was obtained by V. Semyanisty
  in \cite[formulas (9) and (11)]{Semyanistyi:1961},
 \begin{align}\label{SEMR:Semyanistyi}
  f_{even}({\bm \theta})=-\frac{ 1}{4\pi} \int_{{\mathbb S}^2}
   \frac{ 1}
      {({\bm \theta} \centerdot{\bm \eta})^2}
      \{ {\mathcal F} f \}(\bm \eta)
   {\rm d}{\bm \eta},
 \end{align}
 where  the ${\rm d}{\bm \eta}$  is  the surface measure
  on $\mathbb S^2$ with normalization   $\int_{\mathbb S^2}{\rm d}{\bm
  \eta}=4\pi$ and
  integral is understood in the regularized sense.

  In  \cite[p. 99]{Helgason:1999} S. Helgason gives
  for (\ref{SEMR:Funk--Minkowski-scalar}) the inversion formula of filtered
  back-projection type
  \begin{align} \label{SEMR:Helgason}
  f_{even}({\bm \theta}) = \frac{1}{2 \pi }
  \frac{d}{du}\int_0^u
  \int_{{\mathbb S}^2}
  \{ {\mathcal F}f \}({\bm \eta})
  \delta \Big ({\bm \eta} \centerdot{\bm \theta} -\sqrt{1-v^2} \Big )
   {\rm d}{\bm \eta} \frac{v d v}{\sqrt{u^2-v^2}}  \Big |_{u=1}
   \, ,
  \end{align}
  where  $\delta$ denotes the the Dirac delta function.

 Another example of inversion formula is due to
 B. Rubin ~\cite{Rubin:2002, Rubin:2003} %{Rubin:2015}
 \begin{align} \label{SEMR:Rubin}
 f_{even}({\bm \theta})=
 \frac{1}{4\pi} \int_{{\mathbb S}^2}\{ {\mathcal F}f \}({\bm \eta})
 {\rm d}{\bm \eta}
 +  \frac{\Delta_{\bm \theta}}{4\pi} \int_{{\mathbb S}^2}
 \ln|{\bm \eta} \centerdot{\bm \theta}|\, \{ {\mathcal F}f \}({\bm \eta}) \,{\rm d}{\bm \eta} \, ,
 \end{align}
 here  ${\Delta}_{{\bm \theta}}$ it  the Laplace--Beltrami operator
(\ref {SEMR:Laplace--Beltrami}).

 In our studies, an important role is played by spherical convolution
 operator ${\mathcal S}$,
 which is the spherical   analogue  of Hilbert transform,
 see  ~\cite{Samko:1983, Samko:2002, Rubin:2014, Kazantsev:2015}.
    \begin{definition}Let  $f \in  C({\mathbb S}^2)$.
    The  spherical convolution operator ${\mathcal S}$ is defined by,
   \begin{eqnarray}\label{SEMR:OperatorS}
      \{ {\mathcal S}v \}({\bm  \theta}) \equiv
     {\mathcal S}_{{\bm  \theta}}v=p.v.
     \frac{1}{4\pi} \int_{{\mathbb S}^2}
   \frac{v({\bm  \eta})}{{{\bm  \theta}} \centerdot{\bm  \eta}}
   \, {\rm d}{\bm  \eta}, \ {{\bm  \theta}}\in {\mathbb S}^2.
   \end{eqnarray}
 \end{definition}
 This transform is odd,
 $\{ {\mathcal S}f\}(-{\bm \theta})=-\{ {\mathcal S}f\}({\bm \theta})$, and
${\mathcal S}$ annihilates all even  functions.

The results of this paper are  formulated below in Theorems
\ref{th:main1} and \ref{th:main2}.
\begin{theorem}\label{th:main1}
For any function $f({\bm \theta}) \in H^{1}({\mathbb S}^2)$ the
following identity  take place
\begin{align}\label{SEMR:representation}
f(\bm \theta)
&=\underbrace{\frac{1}{4\pi}  \int_{{\mathbb S}^2}
  \{ {\mathcal F}f \}(\bm \eta) {\rm d}{\bm \eta}}_{=f_{00}}+
  p.v.
       \frac{1}{4\pi} \int_{{\mathbb S}^2}
   \frac{
    ({\bm \eta} + {\bm \theta}) \, \centerdot
    \Big \{
    \Big [ {\mathcal F}, \nabla   \Big ] f  \Big   \}   (\bm \eta)
      }
   {{\bm \eta} \centerdot{\bm \theta}}\,
   {\rm d}{\bm \eta}
   \\ \nonumber
   &=f_{00} +
 {\mathcal S}_{\bm \theta}
 ({\bm \eta} + {\bm \theta} ) \centerdot \Big [ {\mathcal F},
\nabla \Big ]_{\bm \eta}f .
\end{align}
 Here operators  ${\mathcal F}$ and  $\nabla$
 %${\mathcal S}$
are  the  Funk--Minkowski transform
(\ref{SEMR:Funk--Minkowski-scalar}) and the surface gradient
(\ref{SEMR:grad}),
%and  the spherical convolution operator (\ref{SEMR:OperatorS}),
respectively. Through the square brackets $ [. , . ]$
we, as usual, denoted the commutator
$ \Big [ {\mathcal F}, \nabla \Big ]f ={\mathcal F} \nabla f
  -\nabla{\mathcal F} f$, where the F--M transform  ${\mathcal F} $ is applied to vector function $\nabla f$ by componentwise.
\end{theorem}

 If  we decompose  identity (\ref {SEMR:representation}) on even and odd parts
  then we can write,
  \begin{align} \nonumber
  f_{even}({\bm \theta})&=
   \frac{1}{4\pi}  \int_{{\mathbb S}^2}
  \{ {\mathcal F}f \}(\bm \eta) {\rm d}{\bm \eta}
   - p.v.
    \frac{1}{4\pi} \int_{{\mathbb S}^2}
   \frac{
   {\bm \theta} \, \centerdot \{ \nabla{\mathcal F} f \}(\bm \eta)}
   {{\bm \theta} \centerdot{\bm \eta}}
   {\rm d}{\bm \eta}
   \\ \label{SEMR:inverseEVEN}
    &=
  \frac{1}{4\pi}  \int_{{\mathbb S}^2}
  \{ {\mathcal F}f \}(\bm \eta) {\rm d}{\bm \eta}
  -{{\bm \theta} \, \centerdot}
    {\mathcal S}_{\bm \theta} \nabla{\mathcal F} f  \, ,
  \\ \label{SEMR:inverseODD}
    f_{odd}({\bm \theta})&= p.v.
       \frac{1}{4\pi} \int_{{\mathbb S}^2}
   \frac{
   {\bm \eta}   \centerdot \{ {\mathcal F} \nabla  f \}(\bm \eta)
      }
   {{\bm \eta} \centerdot{\bm \theta}}\,
   {\rm d}{\bm \eta} = {\mathcal S}_{\bm \theta} {\bm \eta}  \centerdot
   {\mathcal F}_{\bm \eta} \nabla f.
  \end{align}

The inversion formulas  for $f_{even}$ and  $f_{odd}$ follow from these equations %equalities
and if  two F--M transformations
 $g({\bm \eta  })= \{{\mathcal F} f \}({\bm \eta})$ and
 ${\bf h}({\bm \eta  }) =\{ {\mathcal F} \nabla f \}({\bm \eta})$
 are known, then  the unknown function $f$   can be  reconstruct
 completely,
\begin{align}\label{SEMR:fful-f}
  %f({\bm \theta})=
  %\underbrace{
  %\frac{1}{4\pi}  \int_{{\mathbb S}^2}
  %g(\bm \eta) {\rm d}{\bm \eta}
  %-{{\bm \theta} \, \centerdot}
  %\Big \{{\mathcal S} \nabla g \Big\}({\bm \theta})
  %  }_{=f_{even}({\bm \theta})}
  %+
  % \underbrace{ \Big \{ {\mathcal S} {\bm \eta}  \centerdot
  % {\bf h}(\bm \eta )
  % \Big\} ({\bm \theta})
  %      }_{=f_{odd}({\bm \theta})} \, .
  f({\bm \theta})=\frac{1}{4\pi}  \int_{{\mathbb S}^2}
  g(\bm \eta) {\rm d}{\bm \eta}
  - p.v.
    \frac{1}{4\pi} \int_{{\mathbb S}^2}
   \frac{
   {\bm \theta} \, \centerdot  \nabla g(\bm \eta)}
   {{\bm \theta} \centerdot{\bm \eta}}\,
   {\rm d}{\bm \eta}+
    p.v.
       \frac{1}{4\pi} \int_{{\mathbb S}^2}
   \frac{
   {\bm \eta}   \centerdot {\bf h}(\bm \eta)
      }
   {{\bm \eta} \centerdot{\bm \theta}}\,
   {\rm d}{\bm \eta} \, .
    \end{align}

The next problem that we will consider
 is   the problem of Helmholtz--Hodge
 decomposition for  a tangential vector field on the sphere ${\mathbb
 S}^2$, see \cite{FreedenSchreiner:2009}.
 The Helmholtz--Hodge decomposition says that we can write any
 vector field tangent to the surface of the sphere as the sum
   of a curl-free component and a divergence-free component
 \begin{align}\label{SEMR:Helmholtz--Hodge}
 {\bf f}({\bm \theta})=\nabla_{\bm \theta}u({\bm \theta})+{\bm \theta} \times
 \nabla_{{\bm \theta}}v({\bm \theta}),
 \end{align}
  where $\nabla_{{\bm  \theta}} $ is the surface gradient on the sphere, and
 rotated gradient
 ${\bm \theta} \times \nabla_{{\bm  \theta}} $
 means the cross-product of the surface gradient of $v$ with the
 unit normal vector ${\bm \theta}$ to the sphere.
  Here   $\nabla_{\bm \theta} {u}$ is called also as inrrotational,
  poloidal, electric or potential field
  and
  $\nabla^{\bot}_{\bm \theta} {v}$ is called as incompressible,
  toroidal, magnetic or stream vector field.
  Scalar functions ${u}$ and  ${v}$ are called velocity potential and stream
  functions, respectively.

 In the next theorem  we show that  decomposition (\ref{SEMR:Helmholtz--Hodge})
 is obtained
 by use of  Funk--Minkowski- transform ${\mathcal F}$ and  spherical convolution
 transform  ${\mathcal S}$.
  \begin{theorem}\label{th:main2}
  Any vector field ${\bf f} \in {\bf L}_{2,tan}({\mathbb S}^2) $ that is
  tangent to the sphere
  can be uniquely decomposed into a sum
 (\ref{SEMR:Helmholtz--Hodge})
  of a surface curl-free component and  a surface
  divergence-free
  component
  with scalar valued functions $u,v \in H^{1}({\mathbb S}^2)/{\mathbb R}$.
 Functions
 $u$ and $v$ are  velocity potential  and stream
 functions  that are calculated unique up to a constant by the formulas
\begin{align}\nonumber
u({\bm \theta})&=\Big [{\mathcal S, {\bm \eta} \,  \centerdot \,  ,
{\mathcal F}} \Big ]_{\bm \theta}  {\bf f}=
 \Big \{  {\mathcal S } {\bm \eta}
\centerdot {\mathcal F} {\bf f}
 \Big \} ({\bm \theta})-
 \Big \{{\mathcal F} {\bm \eta} \centerdot {\mathcal S}{\bf f}
 \Big \}({\bm \theta})
  \\ \label{SEMR:u}
  &=
   {\mathcal S}_{\bm \theta}
    {\bm \eta}
\centerdot {\mathcal F}_{\bm \eta} {\bf f}
 -
 {\mathcal F}_{\bm \theta} {\bm \eta} \centerdot {\mathcal S}_{\bm \eta}{\bf f}
  ,
\\ \nonumber
v({\bm \theta})&= {\bm \theta} \centerdot
 \Big [{\mathcal S, {\bm \eta} \times , {\mathcal F}} \Big ]_{\bm \theta}{\bf  f}
  ={\bm \theta} \centerdot
 \Big \{
 {\mathcal S } {\bm \eta} \times {\mathcal F} {\bf f}
 \Big \}({\bm \theta}) - {\bm \theta} \centerdot
   \Big \{ {\mathcal F}{ \bm \eta} \times {\mathcal S} {\bf f }
   \Big \} ({\bm \theta})
   \\\label{SEMR:v}
   &={\bm \theta} \centerdot
  {\mathcal S }_{\bm \theta} {\bm \eta} \times {\mathcal F}_{ \bm \eta} {\bf f}
  - {\bm \theta} \centerdot
   {\mathcal F}_{\bm \theta}{ \bm \eta} \times {\mathcal S}_{ \bm \eta} {\bf f }
   ,
\end{align}
 where  through $[{\mathcal A},{\mathcal B},{\mathcal C}]$ we denote the generalized commutator,
\[
 [{\mathcal A},{\mathcal B},{\mathcal C}]={\mathcal A}{\mathcal B}{\mathcal C}
  -{\mathcal C}{\mathcal B}{\mathcal A}.
  \]
\end{theorem}
As a consequence of this  theorem, we can obtain formulas for solving two important problems
on the sphere ${\mathbb S}^2$:
$\nabla u={\bf f}$ and $\nabla^{\bot} v={\bf g}$ .
Answers  to solve these problems are
\[
u({\bm \theta})=({\mathcal S }_{\bm \theta} {\bm \eta} \centerdot
  {\mathcal F}_{\bm \eta} -
  {\mathcal F}_{\bm \theta} {\bm \eta} \centerdot {\mathcal S}_{\bm \eta}
     ){\bf f} \  \ \text{for} \ \   \nabla u={\bf f} \in {\bf L}_{2,tan}({\mathbb S}^2)
\]
and
\[v({\bm \theta})={\bm \theta} \centerdot(
  {\mathcal S }_{\bm \theta} {\bm \eta} \times {\mathcal F}_{\bm \eta}
   - {\mathcal F}_{\bm \theta} {\bm \eta} \times {\mathcal S}_{\bm \eta} ){\bf g}
    \  \ \text{for} \  \ \nabla^{\bot} v={\bf g}\in {\bf L}_{2,tan}({\mathbb S}^2).
\]

\section{Basic methods and tools}
%\section{Definitions and Auxiliary Notations Preliminaries}
\subsection{Spherical harmonics (SHs)}
  %%%%%%%%%%%%%%%%%%%%%%%%%%%%%%%%%%%%%%%%%%%%%%%%%%%%%%%%%%%%%%%%%%%%%%%%%%%%%%%%
  In this section  we state some properties of  complex spherical harmonics.
  A spherical harmonic  $Y_{N\ell}$ of degree  $N$ on ${\mathbb S}^2$
  is the
  restriction to ${\mathbb S}^2$ of a homogeneous harmonic polynomial
  of degree $N$ in ${\mathbb R}^3$.

  The Legendre polynomials of the first kind  $P_{N}$ of degree $N\in {\mathbb N}_0$
  or simply Legendre polynomials   are given by the Rodrigues formula
  \[
  P_{N}(t)=\frac{1}{N!2^N}\frac{{\rm d}^N}{{\rm  d}t^N}(t^2-1)^N.
  \]
  We recall that  Legendre polynomials of the first kind $P_N(t)$
  are the orthogonal polynomials  on $(-1,1)$
  with weight function $w(t)=1.$
  We define with  $C_N^{(3/2)}$ the Gegenbauer polynomial of degree $N$ with  parameter
  ${\lambda}=3/2$,
  \[
  C_N^{(3/2)}(t)=\frac{d}{dt}P_{N+1}(t).
  \]

 The following formulas will be used in our calculations
 (\cite{AbramowitzStegun:1972})
  \begin{align}\label{SEMR:P0}
  P_{2j}(0)=(-1)^j\frac{\Gamma(j+1/2)}{\sqrt{\pi} j!}
  =\frac{(-1)^j(2j-1)!!}{(2j)!!},
 \end{align}
 \begin{align}\label{SEMR:PN}
 (N+1)P_{N+1}(0)=-N P_{N-1}(0),
  \end{align}
  \begin{align}\label{SEMR:P0C0}
      C^{(3/2)}_{2j}(0)=\frac{(-1)^j(2j+1)!!}{(2j)!!}  \  \text{or}  \
      C^{(3/2)}_{N-1}(0)={NP_{N-1}(0)}, \ N=2j+1 .
  \end{align}
The following usefull asymptotics holds as $j$ goes to infinity
  \begin{align}\label{SEMR:asymptotic}
    P_{2j}(0) \sim \frac{1}{\sqrt{2j+1}} \ \text{and} \
     \frac{1}{C^{(3/2)}_{2j}(0)}=\frac{1}{(2j+1)P_{2j}(0)}
     \sim \frac{1}{\sqrt{2j+1}} \ \text{if } \
    j \to \infty .
        \end{align}

The associated Legendre functions of the first kind $P^{\ell}_{N}$  for non negative $\ell\geq  0$  are defined as
  \[
  P^{\ell}_{N}(t)=(1-t^2)^{\frac{\ell}{2}}\frac{{\rm d}^\ell}{{\rm  d}t^\ell}P_N(t),
  \]
  where $N,\ell \in {\mathbb N}_0$ with $\ell \leq N$
  and for the  negative order $-\ell,$ $P^{-\ell}_{N}$ are given by
  \[
  P^{-\ell}_{N}(t)=
  (-1)^{\ell}\frac{(N-\ell)!}{(N+\ell)!}P^{\ell}_{N}(t), \  \ell \geq
  0 \, .
  \]
  When the order $\ell=0$, the associated Legendre function becomes a
  polynomial in $t$ and instead being written $P^{0}_N(t)$ it is designated
  $P_N(t)$, the Legendre polynomial.
  %%%%%%%%%%%%%%%%%%%%%%%%%%%%%%%%%%%%%%%%%%%%%%%%%
  The complex SHs $Y_{N\ell}$
  are related to  the associated Legendre functions as follows
  \[
  Y_{N\ell}({\bm \xi})=(-1)^\ell N_{N\ell}{e}^{{\rm i}\ell\varphi}P^{\ell}_{N}(\cos \theta), \  |\ell| \leq  N,
  \]
  where  $N_{N\ell}$ is a normalization constant
  \[
  N_{N\ell}=\sqrt{\frac{2N+1}{4\pi}\frac{(N-\ell)!}{(N+\ell)!}}
  \]
  and the extra factor $(-1)^\ell$ is called the Condon--Shortley phase.

  The $Y_{N\ell}$ are complex-valued polynomials  of the sines and cosines
  of $\theta$ and $\varphi$
  and for complex conjugate functions  the following formula  fulfil
  \[
  \overline {Y_{N \ell}({\bm \xi})}=(-1)^{\ell} Y_{N,-\ell}({\bm \xi}).
  \]
  %Her  $\overline{\centerdot}$ denote the complex conjugate.
  The  parity rule for spherical harmonic is
  \[
  Y_{N\ell}(-{\bm \xi})=(-1)^NY_{N\ell}({\bm \xi}).
  \]
  It is known that the subspace of all spherical  harmonics of degree $N$,
    ${\rm span} \{ Y_{N\ell} \}_{\ell}^{N}$,
    is the eigenspace of the Laplace--Beltrami
  operator %$\triangle_{\bm \xi}$
  (\ref{SEMR:Laplace--Beltrami})
   corresponding to the eigenvalue $-\lambda^2_N=-N(N+1),$
  \[
  \Delta_{\bm \xi}Y_{N\ell}({\bm \xi})=-N(N+1)Y_{N \ell}({\bm \xi}).
  \]
  The dimension of this subspace being $2N+1$,   so one may choose for it an orthonormal basis  in different ways.

 The collection of all  spherical harmonics
 $\Big \{ Y_{N\ell}, |\ell| \leq N  \Big \}_{N=0}^{\infty}$
 forms an orthonormal basis
 for ${L}_2(\mathbb S^2; {\mathbb C})$
 \begin{equation}
 \label{SEMR:orthSH}
 ({Y}_{N_1\ell_1},  Y_{N_2\ell_2})_{L_2(\mathbb S^2)}=
 \int_{{\mathbb S}^2}{Y}_{N_1\ell_1}({\bm \xi})\overline{Y_{N_2\ell_2}({\bm
 \xi})}
     \, {\rm d} {\bm \xi}=
   \delta^{N_1}_{N_2}\delta^{\ell_1}_{\ell_2},
 \end{equation}
 where  $\delta^{i}_{j}$ is the Kronecker symbol and the space  ${L}_2(\mathbb S^2) \equiv {L}_2(\mathbb S^2; {\mathbb C})$ is a
 Hilbert space of square-integrable functions on ${\mathbb S}^2$
 with the hermitian inner product and the  finite norm,
 \begin{equation*}
 %(u,v)  \equiv
  (u,v)_{L_2(\mathbb S^2)}=
  \int_{\mathbb S^2}  {u}({\bm \xi})  \overline{v({\bm \xi})} \,
  {\rm d}{\bm \xi},
 \
 %||u||^2\equiv
 ||u||^2_{L_2(\mathbb S^2)}= (u,u)_{L_2(\mathbb S^2)}.
 \nonumber
 \end{equation*}

 The Fourier coefficients for $u \in L_2({\mathbb S}^2)$ are
 $ u_{N\ell}=(u,Y_{N\ell})_{L_2}$.
 Then, every  function $u \in {L}_2(\mathbb S^2)$ admits a
 spherical harmonics series expansion  in $L_2$--sense
  \begin{align} \label{SEMR:L2}
  u({\bm \xi}) &=\sum^{\infty}_{N=0} \sum_{\ell} u_{N\ell}
  Y_{N\ell}({\bm \xi}),
  \\
  \label{SEMR:L2norm}
 ||u||^2_{L_2(\mathbb S^2)}
  &=\sum^{\infty}_{N=0} \sum_{\ell} |u_{N\ell}|^2  \, .
   \end{align}

   We close this section with Funk--Hecke formula.
  %,  is contained in the next theorem.
  It was first published
  by Funk (1916) and a little later by Hecke (1918).

  \begin{theorem}\label{th:Funk}[The Funk--Hecke Theorem]
  Suppose $f(t)$ $\in L_1(-1,1)$ is an integrable function.
  Then for
  every spherical harmonics of degree $N$ we have
  \begin{equation}
  \label{SEMR:Funk-Heckeformula}
  \int_{{\mathbb S}^2}f({\bm \xi}\centerdot{\bm \eta})Y_{N \ell}({\bm \xi})
  \, {\rm d}{\bm \xi}= 2\pi Y_{N \ell}({\bm \eta})
  \int_{-1}^{1}f(t)P_N(t)\, {\rm d}t,
  \end{equation}
  where  ${\bm \xi}\centerdot{\bm \eta}$ denotes the inner product
  of  unit vectors ${\bm \xi}$ and ${\bm \eta}$,
  $P_N$ denotes the
  $N$th order Legendre polynomial.
   \end{theorem}

   The Funk--Hecke formula is useful in simplifying calculations of certain integrals
  over ${\mathbb S}^2$ and plays an important role in the theory of  spherical harmonics.
  For more details on the Funk--Hecke formula see ~\cite{AtkinsonHan:2012, Samko:2002},
  for example.
  A general overview on spherical harmonics and the  relevant
  problems can be found in the monographs
  ~\cite{AbramowitzStegun:1972, AtkinsonHan:2012, Muller:1980, VarMosKher:1988,
  Muller:1980, DaiXu:2013, FreedenSchreiner:2009, FreedenGutting:2013}.
   %\subsection{Differential operations and vector analysis on the sphere ${\mathbb S}^2$}
%%%%%%%%%%%%%%%%%%%%%%%%%%%%%%%%%%%%%%
\subsection{Surface differential operators on the sphere ${\mathbb S}^2$}
  Here  we briefly recall the definitions and some properties of surface
  differential operators.

  The space  ${\bf L}_2(\mathbb S^2)\equiv {{\bf L}}_2(\mathbb S^2; {\mathbb C})$
 is a Hilbert space of square-integrable vector functions on
 ${\mathbb S}^2$
 with the inner product and the  finite norm,
 \begin{equation*}
 ({\bf u},{\bf v})_{{\bf L}_2(\mathbb S^2)}=
 \int_{\mathbb S^2}  {{\bf u}}({\bm \xi})  \centerdot  \overline{{\bf v}({\bm \xi})}
 \, {\rm d}{\bm \xi},
 \
 ||{\bf u}||^2_{{\bf L}_2(\mathbb S^2)}= ({\bf u},{\bf u})_{{\bf L}_2(\mathbb S^2)}.
 \nonumber
 \end{equation*}

  %%%%%%%%%%%%%%%%%%%%%%%%%%%%%%%%%%%%%%%%%%%%%%%%%%%%%%%%%%%%%%%%%%%%%%%%%%% \begin{definition}
  \begin{definition}
   {\it The tangential gradient or the surface gradient},
  denoted   by $\nabla \equiv \nabla_{\bm \xi}$
  and
  {\it  the tangential rotated gradient (the surface curl-gradient)}, denoted by
  $\nabla^{\bot} \equiv \nabla^{\bot}_{\bm \xi}$,
  are  defined  accordingly  as
  \begin{align}\label{SEMR:grad}
  \nabla_{\bm \xi}u &=
  \frac{\partial u}{\partial \theta} {\bf e}_1(\bm \xi) +
  \frac{1}{\sin \theta} \frac{\partial u}{\partial \varphi} {\bf  e}_2(\bm \xi),
  \\
  \label{SEMR:grad-bot}
  \nabla^{\bot}_{\bm \xi}u &={\bm \xi} \times \nabla_{\bm \xi}u=
 -\frac{1}{\sin \theta} \frac{\partial u}{\partial \varphi} {\bf e}_{1}(\bm \xi)+
   \frac{\partial u}{\partial \theta} {\bf e}_{2}(\bm \xi),
 \end{align}
 where
  ${\bm \xi} ={\bf  i}\sin{ \,\theta}\cos{ \,\varphi}+
  {\bf  j}\sin{ \,\theta}\sin{ \,\varphi}+
  {\bf  k}\cos{ \,\theta}$.
 \end{definition}

 Obviously, we have  ${\bm \xi} \centerdot \nabla_{\bm \xi}u({\bm
 \xi})=0$, ${\bm \xi} \centerdot \nabla^{\bot}_{\bm \xi}u({\bm
 \xi})=0$ and $\nabla u \centerdot \nabla^{\bot}u =0 $, thus $\nabla u$ and  $\nabla^{\bot} u$
 are will be  tangential vector fields on the sphere ${\mathbb S}^2$
 with $\nabla^{\bot}$
 is rotation by $\pi /2 $ in the tangent plane.

 We  must note here that integration by parts formulas on the
 sphere for  operators (\ref{SEMR:grad}) and  (\ref{SEMR:grad-bot}) are differ.
 Namely, for $u, v \in C^1({\mathbb S}^2),$ we have
  \begin{align}\label{SEMR:intgrad1}
     \int_{{\mathbb S}^2}
     u({\bm \xi})\nabla_{\bm \xi}v(\bm \xi)\,
    {\rm d}{\bm \xi}
  &=-\int_{{\mathbb S}^2}
     v({\bm \xi})\nabla_{\bm \xi}u(\bm \xi)
    {\rm d}{\bm \xi}+2  \int_{{\mathbb S}^2}{\bm \xi}
     u({\bm \xi})v(\bm \xi)\,
    {\rm d}{\bm \xi},
  \\
  \label{SEMR:intgradbot}
     \int_{{\mathbb S}^2}
     u({\bm \xi})\nabla^{\bot}_{\bm \xi}v(\bm \xi)\,
    {\rm d}{\bm \xi}
  &=-\int_{{\mathbb S}^2}
     v({\bm \xi})\nabla^{\bot}_{\bm \xi}u(\bm \xi)\,
    {\rm d}{\bm \xi}.
 \end{align}
 \begin{definition}
  In canonical coordinates, the surface divergence
  ${\rm div}_{\bm \xi}$
  of   vector-valued function
  ${\bf v}({\bm \xi})=v^1{\bf  e}_1({\bm \xi})+v^2{\bf  e}_2({\bm \xi})+v^3{\bm \xi}$
  on the sphere ${\mathbb S}^2$  is written as,
 \begin{align}\label{SEMR:div}
 {\rm div}_{\bm \xi}{\bf v}& =\frac{1}{\sin \theta}\left (
 \frac{\partial}{\partial \theta}(v^{1}\sin \theta )+
 \frac{\partial}{\partial \varphi} v^{2} \right )+2v^3 \, .
  \end{align}
For  tangent vector field ${\bf v}$  we define   the scalar surface rotation (or scalar curl operator) ${\rm curl}_{\bm \xi}$
  by
 \begin{align}\label{SEMR:curl}
 {\rm curl}_{\bm \xi}{\bf v}& =-{\rm div}_{\bm \xi}({\bm \xi }\times {\bf v})=
  \frac{1}{\sin \theta}\left (
 \frac{\partial}{\partial \theta}(v^{2}\sin \theta )-
 \frac{\partial}{\partial \varphi} v^{1} \right ) .
  \end{align}
\end{definition}

 If $u \in C^1({\mathbb S}^2)$ and tangential vector field
   ${\bf v} \in {\bf C}^1({\mathbb S}^2)$, then  we have
   integral formulas, which are also understood as inner  products
  \begin{align}
    \label{SEMR:intdiv}
     \int_{{\mathbb S}^2}
     {\bf v}({\bm \xi}) \centerdot \nabla_{\bm \xi}u(\bm \xi)
    \,{\rm d}{\bm \xi}
  &=-\int_{{\mathbb S}^2}
     u({\bm \xi}){\rm div}_{\bm \xi}{\bf v}(\bm \xi)
    \,{\rm d}{\bm \xi}
    \\ \text{or} \ \
    ({\bf v},\nabla u)_{{\bf L}_2({\mathbb S}^2)}&=
    -({\rm div}\,{\bf v},u)_{L_2({\mathbb S}^2)},
\\    \label{SEMR:intcurl}
     \int_{{\mathbb S}^2}
     {\bf v}({\bm \xi}) \centerdot \nabla^{\bot}_{\bm \xi}u(\bm \xi)
    \,{\rm d}{\bm \xi}
  &=-\int_{{\mathbb S}^2}
     u({\bm \xi}){\rm curl}_{\bm \xi}{\bf v}(\bm \xi)
    \,{\rm d}{\bm \xi}
    \\ \text{or} \ \
   ({\bf v},\nabla^{\bot} u)_{{\bf L}_2({\mathbb S}^2)}&=
    -({\rm curl}\,{\bf v},u)_{L_2({\mathbb S}^2)}.
\end{align}

 \begin{definition}
 Finally, we define  the Beltrami operator, which is also  called the
 Laplace--Beltrami operator
 ${\Delta}\equiv {\Delta}_{\bm \xi}$ as
 \begin{align}\label{SEMR:Laplace--Beltrami}
 \Delta_{\bm \xi}u({\bm \xi})={\rm div}_{\bm \xi}\nabla_{\bm \xi}u({\bm
 \xi})\, ,
 \end{align}
 i.e. the divergence of a gradient is the Laplacian.
 \end{definition}
 One easily checks that
%It is easy to verify that
\begin{align}\label{SEMR:Laplace--Beltrami2}
 \Delta_{\bm \xi}u({\bm \xi})={\rm curl}_{\bm \xi}\nabla^{\bot}_{\bm \xi}u({\bm
 \xi})
 \end{align}
 and also
\begin{align*}
  {\rm curl}_{\bm \xi}\nabla_{\bm \xi}u({\bm  \xi})=0,
\ \ {\rm div}_{\bm \xi}\nabla^{\bot}_{\bm \xi}u({\bm  \xi})=0,
 \end{align*}
 thus we say that $\nabla_{{\bm \xi}} u$ is the curl-free,
  but $ \nabla^{\bot}_{{\bm \xi}}u$ is the divergence-free vector
 fields.

 The next formula is Green--Beltrami identity
 or Green's first surface identity, see \cite[Proposition 3.3]{AtkinsonHan:2012},
 \cite[Theorem 4.12]{Michel:2013}: for
 any $u \in C^1({\mathbb S}^2)$ and any $v \in C^2({\mathbb S}^2)$ we have,
 \begin{align}\label{SEMR:L-B}
     \int_{{\mathbb S}^2}
     \nabla_{\bm \xi}{u}({\bm \xi}) \centerdot \nabla_{\bm \xi}v(\bm \xi)
    \,{\rm d}{\bm \xi}
  &=-\int_{{\mathbb S}^2}
     u({\bm \xi}){\Delta}_{\bm \xi}{v}(\bm \xi)\,
    {\rm d}{\bm \xi}
    \\
    \text{or} \ \
     ( \nabla u,\nabla v)_{{\bf L}_2(\mathbb  S^2)}
    &= -( u, \Delta v)_{L_2(\mathbb  S^2)}
     \, .
        \end{align}
 For example, if we take $u=Y_{N_1\ell_1}$  and  $v=Y_{N_2\ell_2}$, then
\begin{align}\label{SEMR:OrtGrad}
     &( \nabla Y_{N_1\ell_1},\nabla Y_{N_2\ell_2})_{{\bf L}_2(\mathbb
     S^2)}
     \\ \nonumber
    &=
     \int_{{\mathbb S}^2}
     \nabla_{\bm \xi}{Y}_{N_1\ell_1}({\bm \xi})
     \centerdot \nabla_{\bm \xi} \overline {Y_{N_2\ell_2}(\bm
     \xi)}\,
    {\rm d}{\bm \xi}
    =-\int_{{\mathbb S}^2}
     Y_{N_1\ell_1}({\bm \xi}){\Delta}_{\bm \xi}\overline{{Y}_{N_2\ell_2}(\bm
     \xi)}
    {\rm d}{\bm \xi}
    \\ \nonumber
&=N_2(N_2+1)\int_{{\mathbb S}^2}
     Y_{N_1\ell_1}({\bm \xi})\overline{{Y}_{N_2\ell_2}(\bm
     \xi)}\,
    {\rm d}{\bm \xi}=N_2(N_2+1)\delta_{N_1}^{N_2}\delta_{\ell_1}^{\ell_2}
    \, .
        \end{align}
%%%%%%%%%%%%%%%%%%%%%%%%%%%%%%%%%%%%%%%%%%%%%%%%%%%%%%%%%%%%%%%%
 For more definitions and properties of
 these differential operators see e.g. ~\cite{ VarMosKher:1988, AtkinsonHan:2012, Nedelec:2001,
 FreedenSchreiner:2009, FreedenGutting:2013}.
%%%%%%%%%%%%%%%%%%%%%%%%%%%%%%%%%%%%%%%%%%%%%%%%%%%%%%%%%%%%%%
\subsection{Two systems of vector spherical harmonics (VSHs)}
There are  vectorial analogues of scalar spherical harmonics called
vector spherical harmonics.  VSHs can be
defined in several ways.
 %In this article we adopt the following definition.
 In this section we give
 definitions and properties of
 the vector spherical harmonics, which are needed  in our work.
 We refer to ~\cite{MorseFeshbach:1953,  VarMosKher:1988, Nedelec:2001,
DerevtsovKazantsev:2009,   FreedenSchreiner:2009} for more details in this theme.

% Details about them can be found in
%For a comprehensive and unified treatment of both scalar and
% vector spherical harmonics we refer to
%~\cite{FreedenSchreiner:2009}. (W. Freeden and M. Schreiner)

\subsubsection{Pure--spin vector spherical harmonics}

Let us now define a complete  orthogonal set of vectors in ${\bf L}_2({\mathbb S}^2)$.
 \begin{definition}
 The vector spherical harmonics (or pure--spin VSHs)
 are arranged in three families:
 ${\bf y}^{(1)}_{N \ell}({\bm \xi}), $
 ${\bf y}^{(2)}_{N \ell}({\bm \xi}) $
 and ${\bf y}^{(3)}_{N \ell}({\bm \xi}).$
 For ${\bm \xi} \in {\mathbb S}^2$ and given a scalar spherical harmonic
 $Y_{N \ell}({\bm \xi})$
 the unnormalized vector spherical harmonics  are the set
 \begin{align}
  {\bf y}^{(1)}_{N \ell}({\bm \xi})&={\bm \xi}Y_{N \ell}({\bm \xi}), \ \
  N \in 0 \cup {\mathbb N},  \label{SEMR:VSH1}
  \\
 %(2)---------------------------------------------------------------------------
 {\bf y}^{(2)}_{N \ell}({\bm \xi})&=\nabla_{{\bm \xi}}Y_{N \ell}({\bm \xi}),
 \ \ N \in {\mathbb N},
 %= \frac{2{\ell}+1}{4\pi}\int\limits_{S^2}
 %\left [{\bm \eta}-({\bm \xi} \centerdot {\bm \eta}){\bm \xi} \right ]
 %Y_{N\ell}({\bm \eta})C^{(3/2)}_{{\ell}-1}({\bm \xi} \centerdot{\bm \eta}){\rm } \, {\rm d} {\bm \eta},
  \label{SEMR:VSH2}
  \\
 %(3)--------------------------------------------------------------------------
  {\bf y}^{(3)}_{N \ell}({\bm \xi})
  &={\bm \xi}\times {\bf y}^{(2)}_{N \ell}({\bm \xi})
   =\nabla^{\bot}_{\bm \xi}Y_{N \ell}({\bm \xi}),
    \ \ N \in {\mathbb N}.
  %=\frac{2{\ell}+1}{4\pi}{\bm \xi}\times\int_{S^2}
  %{\bf y}^{(1)}_{N\ell}({\bm \eta})
  %C^{(3/2)}_{{\ell}-1}({\bm \xi} \centerdot{\bm \eta}){\rm } \, {\rm d}{\bm \eta}.
  \label{SEMR:VSH3}
 \end{align}
\end{definition}
The pure--spin VSHs form a complete set of orthogonal vector functions on the
surface
 of a sphere ${{\mathbb S}^2}$ with the inner product of the ${\bf L}_2({{\mathbb S}^2})$ space, see \cite [Theorem 5.2.7]{FreedenGutting:2013}.
 % and may orgenized be rearranged as so called
 % "pure-sin vector spherical harmonics" and
 % "pure-orbit vector spherical harmonics".
 %Sometimes we need to   normalize VSHs (\ref{SEMR:VSH1})--(\ref{SEMR:VSH3}).
 Clearly, $||{\bf y}^{(1)}_{N\ell}||_{{\bf L}_2({\mathcal S}^2)}=1.$
  To calculate the  norms  of vector functions ${\bf y}^{(2)}_{N \ell}$ and   ${\bf y}^{(3)}_{N \ell}$, we can use  (\ref{SEMR:OrtGrad}). Therefore, the normalizing vector  harmonics
  or orthonormal  system of  VSHs are
 \begin{align*}{\bf y}^{(1)}_{N\ell}, \
 \widetilde{\bf y}^{(2)}_{N \ell}={\bf y}^{(2)}_{N \ell}/
\sqrt{N(N+1)}\, , \
 \widetilde{\bf y}^{(3)}_{N \ell}={\bf y}^{(3)}_{N
 \ell}/\sqrt{N(N+1)}\, .
 \end{align*}
 Each vector function ${\bf f} \in {\bf L}_2({\mathbb S}^2)$ has
 the Fourier expansion
 \begin{align*}
   {\bf f }({\bm \xi})&=f_{1,00} {\bf y}^{(1)}_{00}({\bm \xi})
   +\sum^{\infty}_{N=1} \sum_{\ell}
  f_{1,N\ell} {\bf y}^{(1)}_{N \ell}({\bm \xi})+
 f_{2,N\ell} \widetilde{\bf y}^{(2)}_{N \ell}({\bm \xi})
 + f_{3,N\ell}\widetilde{\bf y}^{(3)}_{N \ell}({\bm \xi}),
 \\
 ||{\bf f }||^2_{{\bf L}_2({\mathbb S}^2)}&=|f_{1,00}|^2+\sum^{\infty}_{N=1}
 \sum_{\ell} |f_{1,N\ell}|^2 +|f_{2,N\ell}|^2 + |f_{3,N\ell}|^2 .
  \end{align*}

  The hermitian inner products are then given by
\begin{align*}
 ({\bf f },  {\bf h })_{{\bf L}_2({\mathbb S}^2)}
 = f_{1,00} \overline { h}_{1,00}+ \sum^{\infty}_{N=1}
  \sum_{\ell}f_{1,N \ell} \overline { h}_{1,N\ell}+ {f_{2,N \ell} \overline { h}_{2,N\ell}}+ f_{3,N\ell}\overline{h}_{3,N\ell} .
  \end{align*}

\subsubsection{Pure--orbit vector spherical harmonics}

An alternative orthogonal basis in the space ${\bf L}_2({\mathbb S}^2)$
 is the system of pure--orbit VSHs
 $\{{\bf h}^{(e)}_{00}, {\bf h}^{(e)}_{N \ell}, {\bf h}^{(i)}_{N \ell}, {\bf y}^{(3)}_{N \ell},
 \ |\ell| \leq N \}_{N=1}^{\infty}$,
where vector functions ${\bf h}^{(e)}_{N \ell}$ and ${\bf h}^{(i)}_{N \ell}$
defined by
\begin{align}
 \label{SEMR:he0}
  {\bf h}^{(e)}_{00}
  &= -{\bf y}^{(1)}_{00},
\\
 \label{SEMR:he}
  {\bf h}^{(e)}_{N \ell}
  &= -(N+1){\bf y}^{(1)}_{ N \ell}
   + {\bf y}^{(2)}_{N\ell}, \ \ N \in {\mathbb N},
   \\
      \label{SEMR:hi}
    {\bf h}^{(i)}_{N \ell}
  &=N{\bf y}^{(1)}_{N \ell}
  + {\bf y}^{(2)}_{N \ell}, \ N \in {\mathbb N}.
 \end{align}

%
%??The set of vectors in (2.27) are the vector-valued spherical harmonics,
%i.e., the traces on ${\mathbb S}^2$ of vectors whose three components are harmonic polynomials.

 The pure-orbit vector spherical harmonics also has a nice  properties,
 in particular, they are  eigenfunctions for the vectorial
 Funk--Minkowski operator ${\mathcal F}$  in the space ${\bf L}_{2,even}({\mathbb S}^2)$ and for
 vectorial Hilbert  operator ${\mathcal S}$ in the space ${\bf L}_{2,odd}({\mathbb S}^2)$,
  see Lemmas \ref{l:F} and \ref{l:S} in the section Proofs.

%%%%%%%%%%%%%%%%%%%%%%%%%%%%%%%%%%%%%%%%%%%%%%%%%
\subsubsection{Tangent vector fields and  Helmholtz--Hodge decomposition}
 Consider the tangent vector field ${\bf f} \in {\bf L}_{2,tan}({\mathbb S}^2)$,
 it  can be written uniquely as
\begin{align*}
   {\bf f }({\bm \theta})&=
    \underbrace{\sum^{\infty}_{N=1} \sum_{\ell}
 f_{2,N\ell} \widetilde{\bf y}^{(2)}_{N \ell}({\bm \theta})}_
 {\text{the curl-free component}}
 +
 \underbrace{\sum^{\infty}_{N=1} \sum_{\ell}
 f_{3,N\ell}\widetilde{\bf y}^{(3)}_{N \ell}({\bm \theta})}_{\text{the divergence-free component}}
    % \sum^{\infty}_{N=1} \sum_{\ell}
 % f_{2,N\ell} \widetilde{\bf y}^{(2)}_{N \ell}({\bm \xi})
 %+ f_{3,N\ell}\widetilde{\bf y}^{(3)}_{N \ell}({\bm \xi})
\\
&=\sum^{\infty}_{N=1} \frac{1}{\sqrt{N(N+1)}}\sum_{\ell}
 f_{2,N\ell}\nabla Y_{N\ell}({\bm \theta})
 + f_{3,N\ell}{\bm \theta} \times \nabla Y_{N\ell}({\bm \theta}).
 %,
 %\\
 %||{\bf f }||^2_{{\bf L}_2({\mathbb S}^2)}&=\sum^{\infty}_{N=1}
 %\sum_{\ell}|f_{2,N\ell}|^2 + |f_{3,N\ell}|^2 .
 \end{align*}

Then formally we have
\begin{align*}
 {\bf f }({\bm \theta})=\nabla\sum^{\infty}_{N=1} \frac{1}{\sqrt{N(N+1)}}\sum_{\ell}
 f_{2,N\ell} {Y}_{N \ell}({\bm \theta})
 +\nabla^{{\bot}}
\sum^{\infty}_{N=1} \frac{1}{\sqrt{N(N+1)}}\sum_{\ell}
 f_{3,N\ell}{Y}_{N \ell}({\bm \theta}),
  \end{align*}
  where according to (\ref{SEMR:Helmholtz--Hodge})
     the velocity
potential and stream functions are
\begin{align*}
u({\bm \theta})&=\sum^{\infty}_{N=1}
\frac{1}{\sqrt{N(N+1)}}\sum_{\ell}
 f_{2,N\ell} {Y}_{N \ell}({\bm \theta}),
\\
v({\bm \theta})&=\sum^{\infty}_{N=1}
\frac{1}{\sqrt{N(N+1)}}\sum_{\ell}
 f_{3,N\ell}{Y}_{N \ell}({\bm \theta}).
 \end{align*}

Another evident approach  consists in solving the Laplace--Beltrami equations
on the sphere
  \begin{align*}
  \Delta_{\bm \theta}{u}(\bm \theta)&={\rm div}_{\bm \theta} \, {\bf f}(\bm \theta),
  \\
  \Delta_{\bm \theta}{v}(\bm \theta)&={\rm curl}_{\bm \theta} \, {\bf f}(\bm \theta).
  \end{align*}
 They  can be solved
  in integral form, for example, involving
 Green�s function with respect to the Laplace--Beltrami
 $\Delta_{\bm \theta}$, see ~\cite[Theorem 4.6.9]{FreedenGutting:2013}.
%%%%%%%%%%%%%%%%%%%%%%%%%%%%%%%%%%%%%%%%%%%%%%%%%%%%%%%%%%%%%%%%%%%%%%
\subsection{Hilbertian Sobolev spaces on the sphere}
%6.2 Spherical Sobolev Spaces ~\cite{Michel:2013}
\subsubsection{ Sobolev scalar functions  on ${\mathbb S}^2$}

The Sobolev space $H^s({\mathbb S}^2)$ %of order $s$
 with a smoothness index
 $s \geq 0$ is defined by (\cite{AtkinsonHan:2012, Michel:2013,
 Rubin:2014, Nedelec:2001, QuellmalzHielscherLouis:2018})
\begin{align*}
H^s({\mathbb S}^2):= \{ u \in {L}_2(\mathbb S^2; {\mathbb C})
  :\sum^{\infty}_{N=0}(1+N(N+1))^s \sum_{\ell} |u_{N\ell}|^2
  < \infty \} .
  \end{align*}
In other words $u \in H^s({\mathbb S}^2)$ if and only if
$(I-\triangle)^{s/2}u \in L_2({\mathbb S}^2)$. The space
$H^s({\mathbb S}^2)$ is a Hilbert space with the hermitian inner
product
\begin{align*}
 (u,v)_{{H}^s({\mathbb S}^2)}
 =\sum^{\infty}_{N=0} (1+N(N+1))^s
 \sum_{\ell} u_{N\ell} \overline{v_{N\ell}}
  \end{align*}
and the induced norm
\[
||u||^2_{H^s({\mathbb S}^2)}
  =\sum^{\infty}_{N=0}(1+N(N+1))^s \sum_{\ell} |u_{N\ell}|^2=
  ||(I-\triangle)^{s/2}u||^2_{L_2({\mathbb S}^2)}\, .
\]
Putting $s=0$ we obtain $H^0({\mathbb S}^2)=L_2({\mathbb S}^2)$.
If $s=1$ then in addition to (\ref{SEMR:L2}), (\ref{SEMR:L2norm})
we have
\[
 \nabla_{{\bm \xi}}u({\bm \xi})=\sum^{\infty}_{N=0}
\sum_{\ell}
 u_{N\ell}\nabla_{{\bm \xi}}
  Y_{N\ell}({\bm \xi}) ,
=\sum_{N=1}^{\infty}\sqrt{N(N+1)}
 \sum_{\ell}(u,Y_{N\ell})_{L_{2}({\mathbb S}^2)}\widetilde{\bf
 y}^{(2)}_{N\ell}\, ,
 \]
 \begin{align*}
 ||\nabla u||^2_{{\bf L}_2({\mathbb S}^2)}=\sum^{\infty}_{N=0}N(N+1) \sum_{\ell}
 |u_{N\ell}|^2  .
  \end{align*}
Thus we can define the Sobolev space $H^1({\mathbb S}^2)$ as (see
 \cite[p. 14]{Nedelec:2001})
\[H^1({\mathbb S}^2)=\{u \in L_{2} ({\mathbb S}^2) :
\nabla u \in {\bf L}_{2} ({\mathbb S}^2) \}
\]
with its  inner product and the  finite  Sobolev norm
\begin{align*}
 (u,v)_{H^1({\mathbb S}^2)}=(u,v)_{L_2({\mathbb S}^2)}+
(\nabla u,\nabla v)_{{\bf L}_2({\mathbb S}^2)},
 \ \
 ||u||^2_{H^1({\mathbb S}^2)}= ||u||^2_{L_2({\mathbb S}^2)}+
  ||\nabla u||^2_{{\bf L}_{2}({\mathbb S}^2)} \, ,
\end{align*}
where $\nabla$ is the surface gradient on the sphere.
 Generally,
if $s=m$ which is a positive integer, we can define the Sobolev
norm via the following formula
($\nabla$-definition of Sobolev spaces)
\[||u||^2_{H^s({\mathbb S}^2)}=(u, v)_{L_2({\mathbb S}^2)}+
\sum_{k=1}^{m}(\nabla^k u, \nabla^k v)_{{\bf L}_2({\mathbb S}^2)}.
\]
If we s  consider  a closed linear subspace $H^1({\mathbb S}^2)/{\mathbb R})
 \subset H^1({\mathbb S}^2)$,
\[
 H^1({\mathbb S}^2)/{\mathbb R}=\{ u\in H^1({\mathbb S}^2) : \int_{{\mathbb S}^2}
  {u}({\bm \xi}) \,{\rm d}{\bm \xi}=0  \} ,
\]
 then due to a Poincar\'{e} inequality for all
 $u \in H^1({\mathbb S}^2)/{\mathbb R} $ we can  define  an equivalent norm
 for $H^1({\mathbb S}^2)/{\mathbb R}$
\[
||u||^2_{H^1({\mathbb S}^2)/{\mathbb R}}=|| \nabla u||_{{\bf L}_{2}({\mathbb S}^2)} \,,
\]
 such that $H^1({\mathbb S}^2)/{\mathbb R}$ becomes a Hilbert space with
 the inner product
\[
(u,v)_{H^1({\mathbb S}^2)/{\mathbb R}}=(\nabla u,\nabla v)_{{\bf L}_{2}({\mathbb S}^2)} \,.
\]

 %Q. T. Le Gia,
 %Numerical solutions of a boundary value problem on
 %the sphere using radial basis functions
%1) If $s > 1 $ then $H^s({\mathbb S}^2)$ is embedded into
% $C({\mathbb S}^2)$.
%It is also known that if $s>1$, the Sobolev embedding theorem
%gives $H^s \subset C({\mathbb S}^2)$.
%Theorem 6.16.
%2) For all $s \in R$ with $s > 2$, $H^{s}({\mathbb S}^2) \subset
%C^{1}({\mathbb S}^2)$. Thus we  have that $\nabla_{\bm \xi} f$
%exists on  ${\mathbb S}^2$ and is continuous.
%Theorem 6.17.
% 3) Every  $f \in H^{s}({\mathbb S}^2)$ for $s > 2$ is
%Lipschitz continuous  on ${\mathbb S}^2$.

For more details on these spaces, we refer the reader to
     ~\cite{AtkinsonHan:2012},
     ~\cite[Theorems 4.12 and 6.12]{Michel:2013},
     ~\cite[p. 41]{Nedelec:2001} \,.
%[36, 38].
   %(see [24, Chapter 1, Remark 7.6])
\subsubsection{Sobolev tangent vector fields on ${\mathbb S}^2$}
 For tangential vector fields
 we have the vectorial Sobolev space ${\bf H}^{s}_{tan}({\mathbb S}^2)$,
 which is the set of all $f\in {\bf  L}_{2,tan}({\mathbb S}^2)$ such that
\[
||{\bf f }||^2_{{\bf H}^s_{tan}({\mathbb S}^2)}=
\sum^{\infty}_{N=1}(1 + N(N + 1))^s
 \sum_{\ell}|f_{2,N\ell}|^2 + |f_{3,N\ell}|^2.
\]

For the scale of Sobolev spaces ${\bf H}_{tan}^{s}({\mathbb S}^2)$
there is a  Helmholtz--Hodge  decomposition (\cite[Theorem 4.1]{Cantor:1981})
\[{\bf H}_{tan}^{s}({\mathbb S}^2)=\nabla \left (H^{s+1}({\mathbb S}^2)/{\mathbb R} \right )
 \oplus {\rm ker}({\rm div})
  ={\bf H}^{s}_{tan,curl}({\mathbb S}^2)+{\bf H}^{s}_{tan,div}({\mathbb S}^2), \ s \geq 0 .
\]
Here we denote by ${\bf H}^{s}_{tan,div}({\mathbb S}^2)$
  and  ${\bf H}^{s}_{tan,curl}({\mathbb S}^2)$
   the divergence-free and curl-free subspaces of ${\bf H}^{s}_{tan}({\mathbb S}^2)$, respectively.

 %can be  split into two categories referred
 %to as the toroidal ${u}$ and spheroidal ${v}$ functions
 %(in the physical literature the former are also called �magnetic�
 %or �stream�, and the latter �poloidal�, �potential�, or
 %�electric�).
 Another words vector field tangent to the sphere
 ${\bf f} \in {\bf H}^{s}_{tan}({\mathbb S}^2)$ can be uniquely decomposed
 into surface curl-free and surface divergence-free
 components:
\[
{\bf f}=\nabla{u}+ \nabla^{\bot}{v}, \
 \int_{{\mathbb S}^2}{u} \, {\rm d}{\bm \xi}=
 \int_{{\mathbb S}^2}{v} \, {\rm d}{\bm \xi}=0,
\]
where functions $u,v \in H^{s+1}({\mathbb S}^2)/{\mathbb R}$.
We can
define its ${\bf H}^s$ norm, among other equivalent versions, as
\[
 ||{\bf f}||^2_{{\bf H}^s({\mathbb S}^2)}=||u||^2_{H^{s+1}({\mathbb S}^2)} +||v||^2_{H^{s+1}({\mathbb S}^2)}\, .
 \]
\subsection{Fourier multiplier  and spherical convolution operators}

\subsubsection{Fourier multiplier operators}
 Here we define Fourier multiplication operators.
\begin{definition}
The operator  $\Lambda: L_2({\mathbb S}^2) \to L_2({\mathbb S}^2)$
is called the Fourier  multiplier operator with corresponding
sequence of multipliers $\{ \lambda_N \}_{N=0}^{\infty}$
 if operator $\Lambda$ acts  on a function $u \in
L_2({\mathbb S}^2)$ by the formula
 \[
 \{\Lambda u \}({\bm \xi})\equiv \Lambda_{{\bm \xi}} u
 = \sum_{N=0}\lambda_{N}\sum_{\ell}u_{N\ell}Y_{N\ell}({\bm \xi}),
 \]
 %is called the operator with multiplier $\lambda_{N}$,
 where $u_{N\ell}$ denote the Fourier coefficients of $u$ with respect to
 the spherical harmonics,
 \[
 u({\bm \xi})=\sum_{N=0}\sum_{\ell}u_{N\ell}Y_{N\ell}({\bm \xi}).
 \]
\end{definition}
 The sequence of multipliers
 % be a sequence such that
 %Fourier-Laplace  multiplier
 $\{ \lambda_N \}_{N=0}^{\infty}$
 gives complete information about properties of
 operator $\Lambda$, especially  the behavior and asymptotics  of multipliers at
 infinity. It is not hard to see that a multiplier operator on
 $L_2({\mathbb S}^2)$ is bounded if and only if its sequence of multipliers is
 bounded.
% As a generalization
%of spherical convolution operators $\Lambda_u$
The works of many authors are  devoted to the study of such
operators, see \cite{Samko:1983, Rubin:2003, DaiXu:2013}.

\subsubsection{Spherical convolution operators}
An important example of the multiplier operator will be a spherical
convolution operator.
\begin{definition}
 The spherical convolution
 $K *u$  of $K \in
 L_2(-1,1)$ with a function $u \in L_2({\mathbb S}^2)$ is
 defined as
 \[
 (K *u) ({\bm \xi})= \int_{{\mathbb S}^2}
 K({\bm \xi} \centerdot{\bm  \eta})u({{\bm  \eta}})
   \, {\rm d}{\bm  \eta}, \, {\bm \xi}\in S^2,
\]
${\rm d}{\bm \eta}$ is the %(unique)
rotation invariant measure, normalized so that
 $\int_{\mathbb S^2}{\rm d}{\bm \eta}=4\pi$ ---
 the surface area of ${\mathbb S}^2$.
  We recall that  ${\bm \eta} \centerdot {\bm \xi}$ is the usual pointwise  inner
  product.
\end{definition}

By  the Funk�-Hecke formula in Theorem \ref{th:Funk} we have
the sequence of multipliers $\{ \lambda_N \}_{N=0}^{\infty}$
\[
 \{K * Y_{N \ell}\}
({\bm \xi})=2{\pi}
 Y_{N \ell}({\bm \xi})
\int_{-1}^{1} K(x) P_{N}(x) \,dx =\lambda_{N} Y_{N \ell}({\bm
\xi}).
\]

%as well as the spherical cosine transform from (1.3).

\subsubsection{Funk's inversion formula for the F--M transform}
 %inversion by decomposition into spherical harmonics special functions

In \cite{Funk:1913} Funk   showed that Funk--Minkowski-
  transform (\ref{SEMR:Funk--Minkowski-scalar})  is the Fourier multiplier
 operator with multiplicators
 $\lambda_{2j}=P_{2j}(0)$,
 \[
 \{ {\mathcal F} Y_{N \ell} \}({\bm \xi})=P_{2j}(0)Y_{2j,\ell}({\bm\xi}),
\]
 and  asymptotics $\lambda_{2j}=P_{2j}(0) \sim (2j+1)^{-1/2}$  if  $j \to \infty$
 (\cite{AbramowitzStegun:1972}).
 Hence any even function $f_{even} \in C^{\infty}({\mathbb S}^2)$
 can be reconstructed explicitly from
 its Funk--Minkowski transform by the formula
 \begin{align*}
  f_{even}({\bm \xi}) =\sum^{\infty}_{j=0} \sum_{\ell} f_{2j,\ell}
  Y_{2j,\ell}({\bm \xi})=
  \sum^{\infty}_{j=0} \sum_{\ell} \frac{( {\mathcal F} f_{even},
  Y_{2j,\ell})_{L_2(S^2)}}{P_{2j}(0)}
  Y_{2j,\ell}({\bm \xi})\, ,
 \end{align*}
where
\begin{align*}
 ( {\mathcal F} f_{even},  Y_{2j,\ell})_{L_2({\mathbb S}^2)}
 =P_{2j}(0)f_{2j,\ell} \, .
  \end{align*}
 The following mapping property of the Funk--Minkowski transform
 between Sobolev spaces was shown by R.\,S. Strichartz in
 ~\cite[Lemma 4.3]{Strichartz:1981}
 : operator
 \[
  {\mathcal F}:H_{even}^s({\mathbb S}^2) \to H_{even}^{s+1/2}({\mathbb S}^2), \ s\geq 0
 \]
  is continuous and bijective, see also ~\cite{GoodeyWeil:1992, QuellmalzHielscherLouis:2018}.
    %Its  nullspace consists of all odd functions  $f({\bm \xi}) = -f(-{\bm \xi})$.
%%%%%%%%%%%%%%%%%%%%%%%%%%%%%%%%%%%%%%%%%%%%%%%%%%%%%%%%%%
\subsubsection{The spherical convolution operator ${\mathcal S}$}
  Now consider the spherical convolution operator ${\mathcal S}$,
  which defined by formula  (\ref{SEMR:OperatorS}), we repeat it
     \begin{eqnarray*}
      \{ {\mathcal S}v \}({\bm  \xi}) \equiv
     {\mathcal S}_{{\bm  \xi}}v=
    \frac{1}{4\pi} \{ x^{-1} *u\} ({\bm \xi})
=   \frac{1}{4\pi} \int_{{\mathbb S}^2}
   \frac{v({\bm  \eta})}{{{\bm  \xi}} \centerdot{\bm  \eta}}
   \, {\rm d}{\bm  \eta}, \ {{\bm  \xi}}\in {\mathbb S}^2.
   \end{eqnarray*}

      The operator ${\mathcal S}$ does not
   exist as an absolutely convergent integral and
   should be understood in the principal value sense,
   see  ~\cite{Samko:2002, Rubin:2014},
   \[
   \{ {\mathcal S}v \}({{\bm  \xi}})=\lim_{\varepsilon \to 0}\frac{1}{4\pi}
   \int_{|{{\bm  \xi} } \centerdot{\bm  \eta}|> \varepsilon}
   \frac{v({\bm  \eta})}
   {{{\bm  \xi} } \centerdot{\bm  \eta}}
   \, {\rm d}{\bm  \eta}=p.v. \frac{1}{4\pi} \int_{{\mathbb S}^2}
   \frac{v({\bm  \eta})}{{{\bm  \xi}} \centerdot{\bm  \eta}}
   \, {\rm d} {\bm  \eta}.
   \]

  The operator ${\mathcal S}$ is  considered as  operator from
  $L_2({\mathbb S}^2)$ into $L_2({\mathbb S}^2)$ and
  can be regarded as the spherical analogue
  of the Hilbert transform, \cite{Rubin:2014}.
    Evidently, that for even spherical harmonics
  $ \{{\mathcal S}Y_{2j,\ell} \}({\bm  \xi})=0$,
  so we can  consider this operator  only
  on  the subspace of odd SHs, $L_{2,odd}({\mathbb S^2})$.
  %%%%%%%%%%%%%%%%%%%%%%%%%%%%%%%%%%%%%%%%%%%%%%%%%%%%%%
 %The Fourier-Laplace multiplier {hm} of H has the form:

 \begin{proposition} [\cite {Rubin:2014, Kazantsev:2015}]
  The spherical analogue
  of the Hilbert transform  (\ref{SEMR:OperatorS})
 \[
 {\mathcal S}: L_{2,odd}({\mathbb S^2}) \to L_{2,odd}({\mathbb S^2})
 \]
 is a compact operator and a multiplier operator on $L_{2,odd}({\mathbb S^2})$ with
 corresponding sequence of Fourier-Laplace multipliers
 $ \Big \{\frac{1}{C^{(3/2)}_{N-1}(0)}=\frac{1}{NP_{N-1}(0)}, \ N=2j+1
 \Big \}_{j=0}^{\infty}$,
 \begin{eqnarray}
  \label{SEMR:SY}
    \{ {\mathcal S}Y_{N\ell} \}({\bm  \xi})
    &=\frac{1}{C^{(3/2)}_{N-1}(0)}Y_{N\ell}({\bm  \xi}) =\frac{1}{NP_{N-1}(0)}Y_{N\ell}({\bm  \xi}), \ N=2j+1 \
  \end{eqnarray}
and asymptotics
\begin{equation}
    \frac{1}{C^{(3/2)}_{2j}(0)}=\frac{1}{(2j+1)P_{2j}(0)}
    \sim \frac{1}{\sqrt{2j+1}} \ \text{if } \
    j \to \infty .
        \end{equation}
   \end{proposition}

 The operator ${\mathcal S}$, as well as the operator
  ${\mathcal F}$,
  \[
  {\mathcal S}: H_{odd}^s({\mathbb S}^2) \to H_{odd}^{s+1/2}({\mathbb
  S}^2), \ s \geq  0
 \]
  is continuous and bijective in the scale of Sobolev spaces
$H_{odd}^s({\mathbb S}^2)$, see  \cite[Proposition 3.2]{Rubin:2014}.

   %For the proof of compactness we  show that
   %the singular values  decays to zero as $j$ tends to infinity.
   %It is easy verify that the sequence $\frac{1}{\left (C^{(3/2)}_{N}(0) \right )^2}$ is
   %strictly monotonically decreasing to 0 but the sequence
   %$\frac{1}{a_k}=\frac{4j+3}{\left (C^{(3/2)}_{N}(0) \right )^2}$
   %is strictly monotonically increasing to $\pi$.

  %So, we have singular values ~\cite{AbramowitzStegun:1972}
  %\[
  %\frac{(-1)^j}{C^{(3/2)}_{2j}(0)}  \sim \frac{1}{\sqrt{2j+1}} \ j \to  \infty,
  %\]
  %and consequently  operator ${\mathcal S}$ is  compact operator.

\subsubsection{Analytic family of  fractional integrals and  Funk--Minkowski transform}
 %(Generalized Funk--Radon transform in ~\cite{Louis:2016})  and  analytic family}
 We can write the F--M operator (\ref {SEMR:Funk--Minkowski-scalar}) in the form
 of spherical convolution operator  as follows
\begin{align*}
 & \{ {\mathcal F}{u} \}(\bm \xi)=\frac{1}{2\pi}
  \int_{0}^{2\pi}u \Big ({\bf e}_{1}({\bm \xi})\cos \omega
 +{\bf e}_{2}({\bm \xi})\sin \omega \Big ) \,{\rm d}{\omega}
\\
&=
  \frac{1}{2\pi}\int_{-1}^{1}\delta(t)
  \int_{0}^{2\pi}u
  \Big({\bf e}_{1}({\bm \xi})\cos \omega \sqrt{1-t^2}
   +{\bf e}_{2}({\bm \xi})\sqrt{1-t^2}\sin \omega \Big )\,
    {\rm d}{\omega}\,{\rm d}{t}
 \\
 &= \frac{1}{2\pi}
  \int_{{\mathbb S}^2}\delta({\bm \xi} \centerdot{\bm \theta})
  u({\bm \theta}) \, {\rm d}{\bm \theta}=
\frac{1}{2\pi} \{\delta *u\} ({\bm \xi})\, ,
   \end{align*}
   where $\delta$ is the Dirac delta function.

The papers \cite{Louis:2016, QuellmalzHielscherLouis:2018}
give a definition
 of   the generalized Funk--Radon transform
 $S^{(j)}$ for $u\in C^{\infty}({\mathbb S}^2)$  by
  \[
   \{S^{(j)}u\}({\bm \xi}) =\frac{1}{2\pi}
  \int_{{\mathbb S}^2}\delta^{(j)}({\bm \xi} \centerdot{\bm \theta})
  u({\bm \theta}) \, {\rm d}{\bm \theta}, \ \
  j \in 0 \cup {\mathbb N}
  .
 \]
 Here use the notation from ~\cite{Louis:2016, QuellmalzHielscherLouis:2018}
  and $\delta^{(j)}$  denotes the $j$-th derivative of the Dirac delta
 function and  operator $S^{(0)}$ is the Funk--Minkowski
 transform ${\mathcal F}$.

 The spherical Hilbert type operator ${\mathcal S}$ in (\ref{SEMR:OperatorS}) as well as operators
 ${\mathcal S}^{(j)}$ are the members of analytic family
 of fractional integrals $\{ {\mathcal C}^{\lambda}, \widetilde{\mathcal C}^{\lambda} \}$ defined by \begin{align}
 \{{\mathcal C}^{\lambda} f \}({\bm \theta})& =
 \frac { \Gamma \left ( -\frac{\lambda}{2} \right )}
  {2\pi\Gamma \left (  \frac{1+\lambda}{2} \right )}
    \int_{{\mathbb S}^2}
  f({\bm \sigma})
  |{\bm \theta} \centerdot  {\bm \sigma} |^{\lambda}\,  d{\bm \sigma},
  \\
  \{\widetilde{\mathcal C}^{\lambda} f \}({\bm \theta}) &=
 \frac { \Gamma \left ( \frac{1-\lambda}{2} \right )}
  {2\pi \Gamma \left (1+ \frac{\lambda}{2} \right )}
  \int_{{\mathbb S}^2}
 f({\bm \sigma}) |{\bm \theta} \centerdot  {\bm \sigma} |^{\lambda}
  {\rm sgn}({\bm \theta} \centerdot  {\bm \sigma})\,  d{\bm \sigma},
   \end{align}
   see ~\cite{Rubin2:1998, Rubin:2014}.
 The operators  ${\mathcal C}^{\lambda}$ and $\widetilde{{\mathcal C}}^{\lambda}$
 are  called the  $\lambda$-cosine transforms
 of $f$ with even and odd kernel, respectively. If $f \in C^{\infty}({\mathbb S}^2)$, they extend analytically to all $\lambda \in {\mathbb C}$
 with the only poles $\lambda =0,2,4,...$ for ${\mathcal C}^{\lambda}$
 and $\lambda=1,3,5,...$ for $\widetilde{C}^{\lambda}$ .

% ~\cite{Semyanistyi:1961, Samko:1983, Samko:2002, Rubin:2002, Rubin:2003, Rubin:2008}.
The limit case $\lambda =-1$ corresponds to the Funk--Minkowski
  transform ${\mathcal F}$ and Hilbert spherical transform ${\mathcal S}$,
  ~\cite[Lemma 3.4.]{Rubin:2014}
  \begin{align*}{\mathcal F} \sim
\{ {\mathcal C}^{-1} f \} ({\bm \theta} )& =\frac{1} {2\sqrt{\pi}}
 \int_{{\mathbb S}^2}
 f({\bm \eta}) \delta({\bm \theta} \centerdot  {\bm \eta} )\, d{\bm \eta},
 \\
 {\mathcal S} \sim
 \{\widetilde{\mathcal C}^{-1} f \}({\bm \theta}) &= \frac {1}{2 \pi^{3/2}}
 \int_{{\mathbb S}^2}\frac{f({\bm \eta})}{{\bm \theta} \centerdot  {\bm \eta} }\,
   d{\bm \eta}.
  \end{align*}
The integral operator in the inverse formula
(\ref{SEMR:Semyanistyi}) by  V. Semyanisty also belongs to this family with
$\lambda=-2$,
\begin{align*}
 \{ {\mathcal C}^{-2} f \}({\bm \theta})& = \frac {-1}{4\pi^{3/2}}
 \int_{{\mathbb S}^2}
 f({\bm \eta}) \frac{1}{({\bm \theta} \centerdot  {\bm \eta})^2 }
   d{\bm \eta}.
  \end{align*}
The corresponding operator $\widetilde{\mathcal C }^{-2}$ for
 ${\mathcal C}^{-2}$ is the generalized
Funk--Radon transform %S^{(j)}$
\begin{align*}S^{(1)} \sim
 \{ \widetilde{\mathcal C }^{-2} f \}({\bm \theta}) &= \frac {-1}{4 \sqrt{\pi}}
 \int_{{\mathbb S}^2}f({\bm \eta})\delta^{\prime}({\bm \theta} \centerdot  {\bm \eta} )\,
 d{\bm \eta}.
 \end{align*}

If for an analytic continuation we use formulas,
  see for example \cite{GelfandShilov:1964},
\begin{align}
 \frac{ |x|^{\lambda}}
 {\Gamma \left (\frac{1+\lambda}{2}\right )}\Big |_{\lambda=-(2m+1)}
 &=\frac{(-1)^m m!}{(2m)!} \delta^{(2m)}(x), \ m=0,1,2,... \ ,
 \\
 \frac{|x|^{\lambda} {\rm sgn}(x)}
 {\Gamma \left (1+ \frac{\lambda}{2} \right)}\Big |_{\lambda=-2m}&=
 \frac{(-1)^m (m-1)!}{(2m-1)!}\delta ^{(2m-1)}(x)
  , \ m=1,2,3,... \ ,
   \end{align}
%\begin{align}
% ({\mathcal C}^{-2m-1} f)({\bm \theta})& =
% \frac { (-1)^m m!\Gamma \left (m+\frac{1}{2} \right )}
%  {2\pi(2m)!}
%    \int_{S^2}
%  f({\bm \sigma})
%  \delta^{(2m)}(
%  {\bm \theta} \centerdot  {\bm \sigma})\,  d{\bm \sigma}, \ m=0,1,2,... \, .
%  \\
%  (\widetilde{\mathcal C}^{-2m} f)({\bm \theta}) &=
% \frac {(-1)^m (m-1)! \Gamma \left ( m+\frac{1}{2} \right )}
%  {2\pi(2m-1)! }
%  \int_{S^2}
% f({\bm \sigma})\delta ^{(2m-1)}({\bm \theta} \centerdot  {\bm \sigma})
%     \,  d{\bm \sigma}, \ m=1,2,3,... \, .
%    \end{align}
   then as the  result,  the following connection between  ${\mathcal S}^{(2m)}$,
       ${\mathcal S}^{(2m+1)}$ and  analytic family $\{ {\mathcal C}^{\lambda}, \tilde{\mathcal C}^{\lambda} \}$
   take place
  \begin{align}{\mathcal S}^{(2m)} \sim
  \{{\mathcal C}^{-2m-1} f \}({\bm \theta})& =
  \frac { (-1)^m \sqrt{\pi}}
  {2\pi 2^{2m}}
    \int_{{\mathbb S}^2}
  f({\bm \sigma})
  \delta^{(2m)}(
  {\bm \theta} \centerdot  {\bm \sigma})\,  d{\bm \sigma}, \ m=0,1,2,... \,  ,
  \\
   {\mathcal S}^{(2m+1)} \sim  \{ \widetilde{\mathcal C}^{-2m} f \}({\bm \theta}) &=
 \frac {(-1)^m \sqrt{\pi}}
  {2\pi2^{2m-1}}
  \int_{{\mathbb S}^2}
 f({\bm \sigma})\delta ^{(2m-1)}({\bm \theta} \centerdot  {\bm \sigma})
    \,  d{\bm \sigma}, \ m=1,2,3,... \, .
   \end{align}
 According to the general theory of
  analytic family $\{
 {\mathcal C}^{\lambda},  \widetilde{\mathcal C}^{\lambda} \}$ on the sphere ${\mathbb S}^2$,
 we can  find  inverse operators of
 ${\mathcal C}^{\lambda}$, $\widetilde{\mathcal C}^{\lambda}$
  by the formulas (see \cite[Proposition 3.1]{Rubin:2014})
\[  {\mathcal C}^{\lambda} {\mathcal C}^{-\lambda-3}f
={\mathcal C}^{-\lambda-3}{\mathcal C}^{\lambda}f=f, \ \text{where} \
\lambda, -\lambda-3 \neq 0,2,4,... \ f\in C^{\infty}_{even}({\mathbb S}^2),
\]
and
\[  \widetilde{\mathcal C}^{\lambda}\widetilde {\mathcal C}^{-\lambda-3}f
=\widetilde{\mathcal C}^{-\lambda-3}\widetilde{\mathcal C}^{\lambda}f=f, \
\text {where} \ \lambda, -\lambda-3 \neq 1,3,5,... \ ,\  f\in C^{\infty}_{odd}({\mathbb S}^2).
\]
%see ~\cite{Semyanistyi:1961, Samko:1983, Samko:2002, Rubin:2002, Rubin:2003, Rubin:2008}.

In the particular case  $\lambda =-1$  we have
${\mathcal F}^{-1} \sim  \left ({\mathcal C}^{-1} \right )^{-1}={\mathcal C}^{-2}$
and it is appropriate to formula (\ref{SEMR:Semyanistyi}) by V. Semyanisty,
see also \cite [Corollary 3.3]{Rubin:2008}. %~\cite{Semyanistyi:1961}.
If we  apply (formally) the integration by parts formula
(\ref{SEMR:intgrad1}) to  (\ref{SEMR:inverseEVEN}), then we  get
\begin{align*}
 \frac{ {\bm \theta} \, \centerdot}{4\pi} \int_{{\mathbb S}^2}
   \frac{
   \{ \nabla{\mathcal F} f \}(\bm \eta)}
   {{\bm \theta} \centerdot{\bm \eta}}
   {\rm d}{\bm \eta}&=-
 \frac{ {\bm \theta} \, \centerdot}{4\pi} \int_{{\mathbb S}^2}
   \nabla \frac{1}
      {{\bm \theta} \centerdot{\bm \eta}}
      \{ {\mathcal F} f \}(\bm \eta)
   {\rm d}{\bm \eta}
   +\frac{ {\bm \theta} \, \centerdot}{2\pi}
   \int_{{\mathbb S}^2}
   \frac{{\bm \eta}
   \{{\mathcal F} f \}(\bm \eta)}
   {{\bm \eta} \centerdot{\bm \theta}}
   {\rm d}{\bm \eta}
\\&=
 \frac{ {\bm \theta} \, \centerdot}{4\pi} \int_{{\mathbb S}^2}
   \frac{{\bm \theta}-({\bm \theta} \centerdot{\bm \eta}){\bm \eta}}
      {({\bm \theta} \centerdot{\bm \eta})^2}
      \{ {\mathcal F} f \}(\bm \eta)
   {\rm d}{\bm \eta}
   +\frac{ 1}{2\pi}
   \int_{{\mathbb S}^2}
     \{ {\mathcal F} f \}(\bm \eta)
   {\rm d}{\bm \eta}
 \\
 &=
 \frac{ 1}{4\pi} \int_{{\mathbb S}^2}
   \frac{ 1-({\bm \theta} \centerdot{\bm \eta})^2}
      {({\bm \theta} \centerdot{\bm \eta})^2}
      \{ {\mathcal F} f \}(\bm \eta)
   {\rm d}{\bm \eta}
   +\frac{ 1}{2\pi}
   \int_{{\mathbb S}^2}
     \{ {\mathcal F} f \}(\bm \eta)
   {\rm d}{\bm \eta}
 \\
&=
 \frac{ 1}{4\pi} \int_{{\mathbb S}^2}
   \frac{ 1}
      {({\bm \theta} \centerdot{\bm \eta})^2}
      \{ {\mathcal F} f \}(\bm \eta)
   {\rm d}{\bm \eta}
   +\frac{ 1}{4\pi}
   \int_{{\mathbb S}^2}
     \{ {\mathcal F} f \}(\bm \eta)
   {\rm d}{\bm \eta}.
   \end{align*}
  Thus, this  formal calculations show that formula (\ref
 {SEMR:inverseEVEN}) corresponds to formula (\ref{SEMR:Semyanistyi}) and
  serves  as  its  regularization
\begin{align*}
-\frac{ 1}{4\pi} \int_{{\mathbb S}^2}
   \frac{ 1}
      {({\bm \theta} \centerdot{\bm \eta})^2}
      \{ {\mathcal F} f \}(\bm \eta)
   {\rm d}{\bm \eta}=
   \frac{1}{4\pi}  \int_{{\mathbb S}^2}
  \{ {\mathcal F}f \}(\bm \eta) {\rm d}{\bm \eta}
   - p.v.
    \frac{1}{4\pi} \int_{{\mathbb S}^2}
   \frac{
   {\bm \theta} \, \centerdot \{ \nabla{\mathcal F} f \}(\bm \eta)}
   {{\bm \theta} \centerdot{\bm \eta}}
   {\rm d}{\bm \eta}   .
\end{align*}
%%%%%%%%%%%%%%%%%%%%%%%%%%%%%%%%%%%%%%%%%%%%%%%%%%%%%%%%%%%%%%%%%%%%%%
\section{Proofs}
In this section we present the proofs of Theorems \ref{th:main1},
 \ref{th:main2}, which will be based on Lemmas
 \ref{l:F} and \ref{l:S}.
In vector case, as in the scalar case,
 the vectorial   Funk--Minkowski transform
 $
 {\mathcal F}: {\bf L}_{2,even}({\mathbb S^2}) \to {\bf L}_{2,even}({\mathbb S^2})
 $
 and
  vectorial    Hilbert type  spherical transform
 $
 {\mathcal S}: {\bf L}_{2,odd}({\mathbb S^2}) \to {\bf L}_{2,odd}({\mathbb S^2})
 $
 are multiplier operators and relevant mapping properties
 between Sobolev spaces are valid.
 The accurate formulations  are given below.

 \begin{lemma}\label{l:F}
  Vectorial   Funk--Minkowski transform
 $
 {\mathcal F}: {\bf L}_{2,even}({\mathbb S^2}) \to {\bf L}_{2,even}({\mathbb S^2})
 $
  is  a multiplier operator
  \begin{align}\label{SEMR:Fhi}
  {\mathcal F}{\bf h}^{(i)}_{N{\ell}}&=P_{N-1}(0) {\bf h}^{(i)}_{N\ell},  \  N=2j+1,
  \\\label{SEMR:Fy}
  {\mathcal F}{\bf y}^{(3)}_{N{\ell}}&=  P_{N}(0){\bf y}^{(3)}_{N\ell}, \ \ N=2j,
  \\\label{SEMR:Fhe}
  {\mathcal F}{\bf h}^{(e)}_{N\ell} &= P_{N+1}(0) {\bf h}^{(e)}_{N\ell}, \ N=2j+1,
  \end{align}
     where ${\bf h}^{(i)}_{N{\ell}}, {\bf y}^{(3)}_{N{\ell}}, {\bf h}^{(e)}_{N\ell}$
  are  pure--orbit vector spherical harmonics
 (\ref{SEMR:he0})--(\ref{SEMR:hi}).
   We have that in the scale of Sobolev spaces operator
   ${\mathcal F}:{\bf H}_{even}^s({\mathbb S}^2) \to {\bf H}_{even}^{s+1/2}({\mathbb S}^2)$, $s\geq 0$ is continuous and bijective.

  If we choose as a basis  pure--spin vector spherical harmonics,
  then following formulas take place
 \begin{align}\label{SEMR:Fy1}
{\mathcal F}{\bf y}^{(1)}_{N\ell}
&=P_{N-1}(0) {\frac{{\bf y}^{(2)}_{N\ell} }{N+1}}\,
 , \ N=2j+1,
  \\\label{SEMR:Fy2}
{\mathcal F}{\bf y}^{(2)}_{N\ell}
 &=  P_{N-1}(0)\Big ( N{\bf y}^{(1)}_{N\ell} +
 \frac{
  {\bf y}^{(2)}_{N\ell} }{N+1} \Big  ), \ N=2j+1
 \, ,
\\\label{SEMR:Fy3}
    {\mathcal F} {\bf y}^{(3)}_{N{\ell}} &=
    P_{N}(0) {\bf   y}^{(3)}_{N\ell} \, , \  N=2j
    \, .
 \end{align}
    \end{lemma}

Similar statements are valid for the operator
${\mathcal S}$
 \begin{lemma}\label{l:S}
 Vectorial spherical convolution   transform
  $
  {\mathcal S}: {\bf L}_{2,odd}({\mathbb S^2}) \to {\bf L}_{2,odd}({\mathbb S^2})
  $
  is  a multiplier operator
  \begin{align}
  \label{SEMR:She0}
  &{\mathcal S} {\bf h}^{(e)}_{00}
  = {\bf h}^{(e)}_{00},
  \\
  \label{SEMR:She}
  &{\mathcal S} {\bf h}^{(e)}_{N \ell}
  =
  \frac{{\bf h}^{(e)}_{N \ell}}{(N+1)P_{N}(0)} ,  \ N=2j,
  \\
  \label{SEMR:Syi}
   &{\mathcal S} {\bf h}^{(i)}_{N \ell}
    =-\frac{1}{(N+1)P_{N}(0)}\frac{N+1}{N}{\bf h}^{(i)}_{N\ell}, \ N=2j,
  \\
  \label{SEMR:Shy3}
  &{\mathcal S}{\bf y}^{(3)}_{N\ell}=\frac{{\bf y}^{(3)}_{N\ell}}{NP_{N-1}(0)}
  , \ N=2j+1  \, ,
  \end{align}
  where ${\bf h}^{(i)}_{N{\ell}}, {\bf y}^{(3)}_{N{\ell}}, {\bf h}^{(e)}_{N\ell}$
  are  pure--orbit vector spherical harmonics
 (\ref{SEMR:he0})--(\ref{SEMR:hi}).
  In the scale of Sobolev spaces operator
  ${\mathcal S}:{\bf H}_{odd}^s({\mathbb S}^2) \to {\bf H}_{odd}^{s+1/2}({\mathbb S}^2),
  \ s\geq 0$  is continuous and bijective.

The images of  pure--spin spherical harmonics under the action of operator $ {\mathcal S}$
are listed below
   \begin{align}\label{SEMR:Sy0}
  {\mathcal S} {\bf y}^{(1)}_{00} &={\bf y}^{(1)}_{00},
    \\
    \label{SEMR:Sy1}
  {\mathcal S} {\bf y}^{(1)}_{N\ell} &=
 \frac{-1}{P_{N}(0)}\frac{{\bf y}^{(2)}_{N \ell}
 }
 {{N(N+1)}},  \ N=2j,
\\ \label{SEMR:Sy2}
  {\mathcal S} {\bf y}^{(2)}_{N\ell}&=
    \frac{-1}{P_{N}(0)} \Big (
     {\bf y}^{(1)}_{N \ell}
    +\frac{ {\bf y}^{(2)}_{N \ell}}  {N(N+1)}
    \Big ), \    N=2j,
    \\ \label{SEMR:Sy3}
  {\mathcal S} {\bf y}^{(3)}_{N\ell}
   &= \frac{1}{P_{N-1}(0)} \frac{{\bf y}^{(3)}_{N\ell}}{N}
  , \ N=2j+1  \, .
    \end{align}
  \end{lemma}

    {\bf Proof} of  Lemma \ref{l:F}. The
     pure--orbit VSHs are expressed through scalar spherical
 harmonics with the help of three term relations, see for example
 ~\cite{DerevtsovKazantsev:2009},
 \begin{align}\label{SEMR:ThreeTerm1}
  &  {\bf h}^{(i)}_{N \ell}=N{\bf y}^{(1)}_{N \ell}({\bm \xi})+{\bf y}^{(2)}_{N \ell}({\bm \xi})
  \\ \nonumber
  &=
  \alpha_1 Y_{N-1,\ell-1}({\bm \xi})
 \left (
 \begin{array}{c}
   1 \\
   {\rm i} \\
   0 \\
 \end{array}
 \right )
 +\beta_1 Y_{N-1, \ell}({\bm \xi})
 \left (
 \begin{array}{c}
   0 \\
   0 \\
   1 \\
 \end{array}
 \right )
 +\gamma_1 Y_{N-1,\ell+1}({\bm \xi})
 \left (
 \begin{array}{r}
    1 \\
   -{\rm i} \\
    0 \\
 \end{array}
 \right ),
\\
\label{SEMR:ThreeTerm1}
 & {\bf h}^{(e)}_{N \ell}=-(N+1){\bf y}^{(1)}_{N \ell}({\bm \xi})+{\bf y}^{(2)}_{N \ell}({\bm \xi})
 \\ \nonumber
 &=
 \alpha_2 Y_{N+1,\ell-1}({\bm \xi})
\left (
 \begin{array}{c}
   1 \\
   {\rm i} \\
   0 \\
 \end{array}
 \right )
 +\beta_2 Y_{N+1, \ell}({\bm \xi})
 \left (
 \begin{array}{c}
   0 \\
   0 \\
   1 \\
 \end{array}
 \right )
 +\gamma_2 Y_{N+1,\ell+1}({\bm \xi})
 \left (
 \begin{array}{r}
    1 \\
   -{\rm i} \\
     0 \\
 \end{array}
 \right ),
\\%3
& {\bf y}^{(3)}_{N \ell}({\bm \xi})=
 \alpha_{3} Y_{N,\ell-1}({\bm \xi})
\left (
 \begin{array}{c}
   1 \\
   {\rm i} \\
   0 \\
 \end{array}
 \right )
 +\beta_{3} Y_{N \ell}({\bm \xi})
 \left (
 \begin{array}{c}
   0 \\
   0 \\
   1 \\
 \end{array}
 \right )
 +\gamma_{3} Y_{N,\ell+1}({\bm \xi})
 \left (
 \begin{array}{r}
    1 \\
   -{\rm i} \\
    0 \\
 \end{array}
 \right ),
\end{align}
 where $\alpha_{i}, \beta_{i}, \gamma_{i} \ (i=1,2,3)$ some
 coefficients. The values of this coefficients  are unimportant
 here, but their accurate expressions can be found in
 ~\cite{DerevtsovKazantsev:2009}.

By applying the operator ${\mathcal F}$ to these three  term
relations  we immediately obtain:
for  $N=2j+1$
\begin{align}\label{SEMR:Fie}
 {\mathcal F} {\bf h}^{(i)}_{N \ell}= P_{N-1}(0){\bf h}^{(i)}_{N \ell},
\ \
{\mathcal F} {\bf h}^{(e)}_{N \ell}=P_{N+1}(0){\bf h}^{(e)}_{N \ell}
\end{align}
and for  $N=2j$
\[
  {\mathcal F}{\bf y}^{(3)}_{N \ell}=P_{N}(0){\bf y}^{(3)}_{N \ell}\, .
\]

Because the multipliers  have asymptotics $P_{2j}(0) \sim \frac{1}{\sqrt{2j+1}}$
 as $j$ goes to infinity,  we
have that ${\mathcal F}:{\bf H}^{s}_{even}({\mathbb S}^2) \to {\bf H}^{s+1/2}_{even}({\mathbb S}^2)$
 is a continuous operator  in the scale of Sobolev spaces,
as in the scalar case.

%from the Sobolev space of even functions
%${\bf H}^{s}_{even}({\mathbb S}^2)$ onto
%the space
%${\bf H}^{s+1/2}_{even}({\mathbb S}^2)$.
%Its inverse is a continuous operator from the space
%${\bf H}^{s}_{even}({\mathbb S}^2)$
%onto the space ${\bf H}^{s-1/2}_{even}({\mathbb S}^2)$

The two equations (\ref{SEMR:Fie})  can  be written  as
\begin{align*}
 \left \{
 \begin{array}{lll}
   N{\mathcal F}{\bf y}^{(1)}_{N \ell}+{\mathcal F}{\bf y}^{(2)}_{N
 \ell}&=
  NP_{N-1}(0){\bf y}^{(1)}_{N \ell}+P_{N-1}(0){\bf y}^{(2)}_{N\ell}
\\
 -(N+1){\mathcal F}{\bf y}^{(1)}_{N \ell}+{\mathcal F}{\bf y}^{(2)}_{N \ell}
 &=-(N+1)P_{N+1}(0){\bf y}^{(1)}_{N \ell}+P_{N+1}(0){\bf y}^{(2)}_{N \ell} \, .
\\
 \end{array}
\right .
\end{align*}

We need to  solve this system  with respect to
 ${\mathcal  F}{\bf y}^{(1)}_{N \ell}$ and
 ${\mathcal F}{\bf y}^{(2)}_{N \ell}$.
 Subtracting the second from the first equation, we obtain
\begin{align*}
&(2N+1){\mathcal F}{\bf y}^{(1)}_{N \ell}=
  (NP_{N-1}(0)+(N+1)P_{N+1}(0)){\bf y}^{(1)}_{N\ell}
  +(P_{N-1}(0)- P_{N+1}(0)) {\bf y}^{(2)}_{N \ell}
\\
&=P_{N-1}(0)
  \left ( N-(N+1)\frac{N}{N+1} \right ){\bf y}^{(1)}_{N\ell}
  +P_{N-1}(0) \left (1   +  \frac{N}{N+1} \right ) {\bf y}^{(2)}_{N \ell}
  \\
  &=P_{N-1}(0)
  \frac{2N+1}{N+1}{\bf y}^{(2)}_{N \ell}.
\end{align*}
 Here  we used the formula (\ref{SEMR:PN}),  $(N+1)P_{N+1}(0)=-N P_{N-1}(0)$,
 thus   we have
\[{\mathcal F}{\bf y}^{(1)}_{N \ell}=P_{N-1}(0)\frac{{\bf y}^{(2)}_{N \ell}}{N+1}, \  \
{\mathcal F}{\bf y}^{(2)}_{N\ell}
 =  P_{N-1}(0)\Big ( N{\bf y}^{(1)}_{N\ell} +
 \frac{
  {\bf y}^{(2)}_{N\ell} }{N+1} \Big ), \ N=2j+1.
\]

{\bf Proof} of  Lemma \ref{l:S}.
By applying  operator ${\mathcal S}$ to three  term relations,
as well as in the previous case, we obtain:
for  $N=2j$
\begin{align*}
    {\mathcal S} {\bf h}^{(i)}_{N \ell}
   & = \frac{{\bf h}^{(i)}_{N \ell}}{(N-1)P_{N-2}(0)}
   =
   -\frac{1}{(N+1)P_{N}(0)}\frac{N+1}{N}{\bf h}^{(i)}_{N\ell},
   \\
   {\mathcal S} {\bf h}^{(e)}_{N \ell}
  &=
  \frac{{\bf h}^{(e)}_{N \ell}}{(N+1)P_{N}(0)}
  \end{align*}
and for  $N=2j+1$
\[
{\mathcal S}{\bf y}^{(3)}_{N \ell}=\frac{{\bf y}^{(3)}_{N \ell}}{NP_{N-1}(0)}\, .
\]
In the first formula we used equality
$(N-1)P_{N-2}(0)=-NP_{N}(0)$.

Continuity of the operator ${\mathcal S}$ in the scale
${\bf H}_{odd}^s({\mathbb S}^2)$
 follows from  asymptotic behavior
$\frac{1}{(2j+1)P_{2j}(0)} \sim \frac{1}{\sqrt{2j+1}}$ if  $j \to \infty$.

The first two equations above are equivalent to the system
\begin{align*}
 \left \{
 \begin{array}{lll}
 N{\mathcal S}{\bf y}^{(1)}_{N \ell}+{\mathcal S}{\bf y}^{(2)}_{N \ell}
  &=
  -\frac{{\bf y}^{(1)}_{N \ell}}{P_{N}(0)}-\frac{{\bf y}^{(2)}_{N \ell}}{N P_{N}(0)}
 \\
  -(N+1){\mathcal S}{\bf y}^{(1)}_{N \ell}+{\mathcal S}{\bf y}^{(2)}_{N \ell}
 &=
  -\frac{{\bf y}^{(1)}_{N \ell}}{P_{N}(0)}+ \frac{{\bf y}^{(2)}_{N  \ell}}{(N+1)P_{N}(0)}\,.
  \\
 \end{array}
\right .
\end{align*}

 Solving this system, we obtain the desired
\begin{align*}
  {\mathcal S}{\bf y}^{(1)}_{N \ell}=
         \frac{-{\bf y}^{(2)}_{N \ell}}{P_{N}(0)N(N+1)} ,
\ \
 {\mathcal S}{\bf y}^{(2)}_{N \ell}=
    \frac{-1}{P_{N}(0)}
   \Big  ( {\bf y}^{(1)}_{N \ell}
       + \frac{{\bf y}^{(2)}_{N \ell}}{N(N+1)} \Big  )  .
\end{align*}
\subsection{Proof Theorem \ref{th:main1}}
%\begin{proof}
We recall some of the basic properties that are implied in our proof.
The Funk--Minkowski transform is even, $\{ {\mathcal F}f\}(-{\bm
\xi})=\{ {\mathcal F}f\}({\bm \xi})$, and ${\mathcal F}
f_{odd} =0$,
 but spherical transform  ${\mathcal S}$ is odd,
   $\{{\mathcal S} f\}(-{\bm \xi})=- \{ {\mathcal S} f\}({\bm
   \xi})$,
and ${\mathcal S} f_{even} =0$. It is obviously that the surface
gradient  $\nabla$, scalar (dot) product ${\bm \eta} \, \centerdot
$ and vector (cross) product ${\bm \eta} \times $ change the
parity.
 We also recall the  parity rules for scalar and vector spherical harmonics :
  $  Y_{N\ell}(-{\bm \xi})=(-1)^NY_{N\ell}({\bm \xi})$,
  ${\bf y}^{(1)}_{N\ell}(-{\bm \xi})= (-1)^{N+1}{\bf y}^{(1)}_{N\ell}({\bm \xi})$,
  ${\bf y}^{(2)}_{N\ell}(-{\bm \xi})= (-1)^{N+1}{\bf y}^{(2)}_{N\ell}({\bm \xi})$,
  ${\bf y}^{(3)}_{N\ell}(-{\bm \xi})= (-1)^{N}{\bf y}^{(3)}_{N\ell}({\bm \xi})$.

 Now we can proceed to our
formula (\ref{SEMR:representation}) and without loss of
generality, we assume that
 $f({\bm \theta}) \in H^{1}({\mathbb S}^2)/{\mathbb R}$, then we have
\begin{align*}
 f({\bm \theta})&=
        \frac{1}{4\pi} \int_{{\mathbb S}^2}
   \frac{
    ({\bm \eta} + {\bm \theta}) \, \centerdot
    \Big \{
    \Big [ {\mathcal F}, \nabla   \Big ] f  \Big   \}   (\bm \eta)
      }
   {{\bm \eta} \centerdot{\bm \theta}}\,
   {\rm d}{\bm \eta}
\\ &=
        \frac{1}{4\pi} \int_{{\mathbb S}^2}
   \frac{
    {\bm \eta}     \centerdot \Big \{  \Big [ {\mathcal F}, \nabla   \Big ]
 f \Big \}
    (\bm   \eta)+
   {\bm \theta} \, \centerdot  \Big \{ \Big [ {\mathcal F}, \nabla   \Big ]
   f \Big \}(\bm \eta)
      }
   {{\bm \eta} \centerdot{\bm \theta}}\,
   {\rm d}{\bm \eta}
   \\
  &=  \frac{1}{4\pi} \int_{{\mathbb S}^2}
   \frac{
    {\bm \eta}   \centerdot  {\mathcal F}_{\bm   \eta} \nabla  f -
\overbrace{ {\bm \eta}   \centerdot  \nabla{\mathcal F}_{\bm   \eta} f }^{=0}
    -
   {\bm \theta}  \centerdot \nabla{\mathcal F}_{\bm   \eta} f
     }
   {{\bm \eta} \centerdot{\bm \theta}}\,
   {\rm d}{\bm \eta}+
       \frac{ {\bm \theta}  \centerdot}{4\pi} \int_{{\mathbb S}^2}
   \frac{
   \overbrace{  {\mathcal F}_{\bm   \eta} \nabla  f}^{eveen}
      }
   {{\bm \eta} \centerdot{\bm \theta}}\,
   {\rm d}{\bm \eta}
      \\
 &=       \frac{1}{4\pi} \int_{{\mathbb S}^2}
   \frac{
    {\bm \eta}   \centerdot {\mathcal F}_{\bm   \eta} \nabla  f-
   {\bm \theta} \centerdot  \nabla {\mathcal F}_{\bm   \eta} f
      }
   {{\bm \eta} \centerdot{\bm \theta}}\,
   {\rm d}{\bm \eta}
   = {\mathcal S}_{\bm \theta} {\bm \eta}   \centerdot {\mathcal F}_{\bm \eta} \nabla  f -
    {\mathcal S}_{\bm \theta}{\bm \theta}  \centerdot \nabla {\mathcal F}_{\bm \eta} f .
   \end{align*}

   It is clear that
   ${\rm ker} ({\mathcal S}_{\bm \theta} {\bm \eta}
    \centerdot {\mathcal F}_{\bm \eta} \nabla)=H^{s}_{even}({\mathbb S }^2)$
   and
     ${\rm ker} ({{\bm \theta}  \centerdot}\,
   {\mathcal S}_{\bm \theta} \nabla{\mathcal F}_{\bm \eta})
 =H^{s}_{odd}({\mathbb S }^2)$.
 From the  Lemmas  \ref{l:F},  \ref{l:S} we have
  ${\mathcal S} {\bm \eta}  \centerdot
   {\mathcal F}_{\bm \eta}  \nabla :H^{s}({\mathbb S }^2)\to H^{s}({\mathbb S }^2)$
   if $s\geq 1$, that looks on the diagram
\[
   \begin{CD} H^{s}({\mathbb S }^2)  @>  \nabla     >>
{\bf H}^{s-1}_{tan}({\mathbb S }^2)
                         @>  {\bm \eta }  \centerdot {\mathcal F}_{\bm \eta}          >>
           H^{s-1/2}({\mathbb S }^2)  @>  {\mathcal S}       >>
             H^{s}({\mathbb S }^2).
   \\
   \end{CD}
\]
Similarly, $ {\bm \theta} \centerdot
  {\mathcal S}_{\bm \theta} \nabla {\mathcal F}:
  H^{s}({\mathbb S }^2)  \to   H^{s}({\mathbb S }^2)$ if $s\geq 1/2$,
  which is also confirmed by the diagram
  \[
  \begin{CD} H^{s}({\mathbb S }^2)  @>  {\mathcal F}      >>
  H^{s+1/2}({\mathbb S }^2)          @>   \nabla        >>
 {\bf H}^{s-1/2}_{tan}({\mathbb S }^2)
           @>  {\bm \theta}  \centerdot {\mathcal S}_{\bm \theta}       >>
     H^{s}({\mathbb S }^2).
   \\
   \end{CD}
 \]
%Since our  proof will consist in a direct verification, for the
%sake of clarity, let's draw a diagram of our calculations.
 For further calculations we take specifically $f={Y}_{N\ell}$, then
 from (\ref{SEMR:Fy2}) and (\ref{SEMR:Fy3}) we have
 \begin{align*} \text{for} \ N=2j &:
 {\mathcal F}_{\bm \eta} \nabla {Y}_{N\ell}= 0,
 \ \
  \nabla {\mathcal F}_{\bm \eta} {Y}_{N\ell}
 =  P_{N}(0) {\bf y}^{(2)}_{N\ell} ({\bm \eta}),
  \\ \text{for} \ N=2j+1 &:
  {\mathcal F}_{\bm \eta} \nabla {Y}_{N\ell}
 =  P_{N-1}(0)\Big ( N{\bf y}^{(1)}_{N\ell}({\bm \eta}) +
 \frac{
  {\bf y}^{(2)}_{N\ell} ({\bm \eta})}{N+1}
 \Big  ), \ \
  \nabla
 {\mathcal F}_{\bm \eta} {Y}_{N\ell}
 =0.
 \end{align*}
Consequently
\begin{align*}
{\bm \eta}   \centerdot  {\mathcal F}_{{\bm \eta}} \nabla  Y_{N\ell} -
   {\bm \theta} \centerdot \nabla{\mathcal F}_{\bm \eta} Y_{N\ell}
 = \left \{
\begin{array}{ll}
   -P_{N}(0)
   {\bm \theta} \centerdot  {\bf y}^{(2)}_{N\ell} ({\bm \eta}), & N=2j \\
   P_{N-1}(0) N{Y}_{N\ell}({\bm \eta}),
  & N=2j+1
  \\
\end{array}
\right.
\end{align*}
and finally,   using formulas (\ref{SEMR:SY}) and
 (\ref{SEMR:Sy2}), we get
\begin{align*}
       \frac{1}{4\pi} \int_{{\mathbb S}^2}
   \frac{
   {\bm \eta}   \centerdot  {\mathcal F}_{{\bm \eta}} \nabla  Y_{N\ell} -
   {\bm \theta} \, \centerdot \nabla{\mathcal F}_{\bm \eta} Y_{N\ell}
         }
   {{\bm \eta} \centerdot{\bm \theta}}\,
   {\rm d}{\bm \eta}=
   \left \{
\begin{array}{ll}
   {Y}_{N\ell}({\bm \theta}), & N=2j \\
   {Y}_{N\ell}({\bm \theta}),
  & N=2j+1 \,  .
  \\
\end{array}
\right.
\end{align*}
So we proved that, if
    $s\geq 1$ then
operator
${\mathcal S}_{\bm \theta} {\bm \eta}   \centerdot {\mathcal F}_{\bm \eta} \nabla -
    {\mathcal S}_{\bm \theta}{\bm \theta}  \centerdot \nabla {\mathcal F}_{\bm \eta}
    :H^{s}({\mathbb S }^2) \to H^{s}({\mathbb S }^2)$
    is identical operator.

%\end{proof}
%%%%%%%%%%%%%%%%%%%%%%%%%%%%%%%%%%%%%%%%%%%%%%%%%%%%
\subsection{Proof Theorem \ref{th:main2}}%the Helmholtz--Hodge
  We have already mentioned two approaches to the solution of
  Helmholtz--Hodge  decomposition problem for
  ${\bf f} \in {\bf L}_{2, tan}({\mathbb S}^2)$.
 Now we proof  formulas (\ref{SEMR:u}) and
(\ref{SEMR:v}) in Theorem  \ref{th:main2}
 for the velocity potential  $u$ and stream functions $v$ of  the Helmholtz--Hodge
  decomposition (\ref{SEMR:Helmholtz--Hodge}),
\begin{align*}
u({\bm \theta})&=
   {\mathcal S }_{\bm \theta} {\bm \eta}
 \centerdot {\mathcal F}_{\bm \eta} {\bf f}
 -
 {\mathcal F}_{\bm \theta} {\bm \eta} \centerdot {\mathcal S}_{\bm \eta}{\bf f},
\\ v({\bm \theta})&= {\bm \theta} \centerdot
 {\mathcal S }_{\bm \theta} {\bm \eta} \times {\mathcal F}_{\bm \eta} {\bf f}
 - {\bm \theta} \centerdot
   {\mathcal F}_{\bm \theta}{ \bm \eta} \times {\mathcal S}_{\bm \eta} {\bf f }.
  \end{align*}

  For the proof it suffices to verify these formulas on the basis elements
 $f={\bf y}^{(2)}_{N\ell}$ and
$f={\bf y}^{(3)}_{N\ell}$. Applying the scalar and cross products
to the the formulas (\ref{SEMR:Fy2})--(\ref{SEMR:Fy3}) and
(\ref{SEMR:Sy2})--(\ref{SEMR:Sy3}), we obtain
\begin{align*}
\text{for} \  N=2j &:
 {\bm \eta} \centerdot  {\mathcal F}_{\bm \eta} {\bf y}^{(3)}_{N{\ell}}= 0,
 \ {\bm \eta} \times {\mathcal F}_{\bm \eta} {\bf y}^{(3)}_{N{\ell}}=
   - P_{N}(0) {\bf   y}^{(2)}_{N\ell}({\bm \eta}),
 \\
 \text{for} \   N=2j+1 &:
 {\bm \eta} \centerdot {\mathcal F}_{\bm \eta}{\bf y}^{(2)}_{N\ell}
 =  P_{N-1}(0) N {Y}_{N\ell}({\bm \eta}), \
 {\bm \eta} \times {\mathcal F}_{\bm \eta}{\bf y}^{(2)}_{N\ell}
 =  P_{N-1}(0)\frac{
  {\bf y}^{(3)}_{N\ell} ({\bm \eta})}{N+1},
       \end{align*}

\begin{align*}
 \text{for} \    N=2j &:
        {\bm \eta} \centerdot
{\mathcal S}_{\bm \eta} {\bf y}^{(2)}_{N\ell}=
    \frac{-1}{P_{N}(0)}{Y}_{N \ell}({\bm \eta}), \ \
      {\bm \eta} \times
 {\mathcal S}_{\bm \eta} {\bf y}^{(2)}_{N\ell}=
    \frac{-1}{P_{N}(0)} \frac{ {\bf y}^{(3)}_{N \ell}({\bm \eta})}  {N(N+1)},
    \\
    \text{for} \  N=2j+1&:
   {\bm \eta} \centerdot {\mathcal S}_{\bm \eta}
   {\bf y}^{(3)}_{N\ell}=0, \
   {\bm \eta} \times {\mathcal S}_{\bm \eta} {\bf y}^{(3)}_{N\ell}
   = \frac{-1}{P_{N-1}(0)} \frac{{\bf y}^{(2)}_{N\ell}({\bm \eta})}{N} .
 \end{align*}

Based on the above, we get
 \begin{align*}
  {\mathcal S }_{\bm \theta} {\bm \eta} \centerdot
  {\mathcal F}_{\bm \eta}  {\bf y}^{(2)}_{N\ell}
  -
  {\mathcal F}_{\bm \theta} {\bm \eta} \centerdot {\mathcal S}_{\bm \eta}
      {\bf y}^{(2)}_{N\ell}
 & = \left \{
 \begin{array}{ll}
  {Y}_{N \ell}({\bm \theta}), & N=2j   \\
  {Y}_{N \ell}({\bm \theta}) , & N=2j+1 ,  \\
 \end{array}
 \right .
 \\
 {\mathcal S }_{\bm \theta} {\bm \eta} \centerdot {\mathcal F}_{\bm \eta} {\bf y}^{(3)}_{N\ell}
 -
 {\mathcal F}_{\bm \theta} {\bm \eta} \centerdot {\mathcal S}_{\bm \eta}{\bf y}^{(3)}_{N\ell}
 &= \left \{
 \begin{array}{ll}
  0, & N=2j \\
  0, & N=2j+1 , \\
 \end{array}
 \right .
 \\
{\bm \theta} \centerdot
   {\mathcal S }_{\bm \theta} {\bm \eta} \times {\mathcal F}_{\bm \eta} {\bf y}^{(2)}_{N\ell}
 -{\bm \theta} \centerdot
 {\mathcal F}_{\bm \theta} {\bm \eta} \times {\mathcal S}_{\bm \eta} {\bf y}^{(2)}_{N\ell}
  &= \left \{
 \begin{array}{ll}
  0, & N=2j \\
  0 & N=2j+1 , \\
 \end{array}
 \right .
 \\
 {\bm \theta} \centerdot
  {\mathcal S }_{\bm \theta} {\bm \eta} \times {\mathcal F}_{\bm \eta} {\bf y}^{(3)}_{N\ell}
 -{\bm \theta} \centerdot
 {\mathcal F}_{\bm \theta} {\bm \eta} \times {\mathcal S}_{\bm \eta}{\bf y}^{(3)}_{N\ell}
  &= \left \{
 \begin{array}{ll}
  {Y}_{N \ell}({\bm \theta}), & N=2j \\
  {Y}_{N \ell}({\bm \theta}), & N=2j+1 . \\
 \end{array}
 \right .
  \end{align*}
\section{Conclusion}
 This paper is devoted to the study of Funk--Minkowski
 transform ${\mathcal F}$
  and  Hilbert type spherical convolution ${\mathcal S}$.  We provide inversion
 formulas for  two F--M transforms ${\mathcal F }f$ and ${\mathcal F} \nabla f$. In this case both  even and
 odd parts of the function $f$ are determined. Also, the formulas
 for decomposition of a tangent vector field on the sphere into
 divergence-free and curl-free parts
 with the participation of operators ${\mathcal F}$ and
 ${\mathcal S}$
 are  derived. In the process
 of obtaining and proving all formulas, the spherical multipliers
 approach  is used.
%Fourier multiplication operators

%A B C D E F G H I K L M N O P Q R S T V X Y Z


\bigskip
\begin{thebibliography}{100}

\bibitem{AbramowitzStegun:1972} M. Abramowitz and I.\,A. Stegun,
{\it Handbook of Mathematical Functions with Formulas, Graphs, and
Mathematical Tables}, 9th printing, New York, Dover, 1972.
%Electronic copy available in
%\url{http://www.math.sfu.ca/~cbm/aands/}
% Abramowitz, M. and Stegun, I. A. (Eds.).
% ``Bernoulli and Euler Polynomials and the Euler-Maclaurin Formula.''
% {\S}23.1 in Handbook of Mathematical Functions with Formulas, Graphs,
% and Mathematical Tables, 9th printing. New York: Dover, pp.~804-806, 1972.

\bibitem{AbouelazDaher:1993}
 A. Abouelaz and R. Daher, {\it Sur la transformation de Radon de
la sphere ${\mathbb S}^d$}, Bull. Soc. math. France,
 {\bf 121}:3 (1993), 353--382.


\bibitem{AtkinsonHan:2012}
 K. Atkinson and W. Han, {\it Spherical Harmonics and Approximations on the Unit
Sphere: An Introduction}, volume 2044 of Lecture Notes in
Mathematics, Springer, Heidelberg, 2012.

%\bibitem{Helgason2008}
%S. Helgason, Geometric Analysis on Symmetric Spaces.  AMS, 2nd
%edition, 2008.

\bibitem{Cantor:1981}
M. Cantor, {\it Elliptic operators and the decomposition of tensor
fields}, Bulletin (new series) of the American mathematical
society,  {\bf 5}:3 (1981), 235--262.

\bibitem{DaiXu:2013}
  F. Dai and Y. Xu, {\it Approximation Theory and Harmonic
 Analysis on Spheres and Balls}, Springer Monographs in Mathematics,
 Springer, New York, 2013.

\bibitem{Daher:2001}
 R. Daher, {\it  Un theoreme de support pour une transformation
de Radon sur la sphere ${\mathbb S}^d$}, C. R. Acad. Sci. Paris,
{\bf 332}:9 (2001), 795--798.

\bibitem{Dann:2010}
S. Dann, {\it On the Funk--Minkowski transform}, (2010),
arXiv:1003.5565.

\bibitem{DerevtsovKazantsev:2009}
E.\, Yu. Derevtsov, S.\,G. Kazantsev, Th. Schuster, {\it
Polynomial bases for subspaces of vector fields in the unit ball.
Method of ridge functions}, J. Inv. Ill-Posed Problems,
  {\bf 15}:1 (2007), 19--55.


\bibitem{Funk:1913}
 P. Funk, {\it \"{U}ber Fl\"{a}chen mit
 lauter geschlossenen geod\"{a}tischen Linien}, Math. Ann.,
 {\bf 74}:2 (1913), 278--300.
%%Mathematische Annalen, 74(2):278--300., 1913.
%(1914) u Natterera?

 \bibitem{Funk:1915}
 P. Funk, {\it Beitrage zur Theorie der Kugelfunktionen}, Math.
Ann., {\bf 77} (1915), 136--152.

\bibitem{Funk:1916}
P. Funk, {\it Uber eine geometrische Anwendung der Abelschen Integralgleichung},
Math. Ann., {\bf 77}  (1916), 129--135.



\bibitem{FreedenSchreiner:2009}
W. Freeden and M. Schreiner, {\it Spherical functions of
mathematical geosciences. A scalar, vectorial, and tensorial
setup}, Berlin: Springer, 2009.


%W. Freeden and M. Schreiner, Spherical Functions of Mathematical Geosciences { A Scalar,
%Vectorial, and Tensorial Setup, Advances in Geophysical and Environmental Mechanics
%and Mathematics, Springer Berlin Heidelberg, 2009.



\bibitem{FreedenGutting:2013}
W. Freeden, M. Gutting, {Special Functions of Mathematical (Geo-)Physics}, Springer Basel, 2013.


%\bibitem{FreedenGervens:1993}
%W. Freeden and T. Gervens, {\it Vector spherical spline
%interpolation basic theory and computational aspects},
% Math. Methods Appl. Sci., {\bf 16} (1993),  151--183.

%%\bibitem{FreedenGervensSchreiner:1998}
%%  W. Freeden, T. Gervens and M.   Schreiner,
%%  {\it Constructive Approximation on the Sphere},
%%  Oxford University Press, Oxford, 1998.

% W. Freeden, T. Gervens, and M. Schreiner, Constructive
%Approximation on the Sphere, in Numer. Math. Sci. Comput., The
%Clarendon Press, Oxford University Press, New York, 1998.

\bibitem{GelfandShilov:1964}
I.\,M. Gelfand and G.\,E. Shilov, {\it Generalized Functions 1,
Properties and Operations}, Academic Press, New York-London, 1964.


\bibitem{GindikinReedsShepp:1994}
S. Gindikin, J. Reeds, and L. Shepp, {\it Spherical tomography and
spherical integral geometry}, In E. T. Quinto, M. Cheney, and P.
Kuchment, editors, Tomography, Impedance Imaging, and Integral
Geometry, volume 30 of Lectures in Appl. Math, pages 83--92.
American Mathematical Society, South Hadley, Massachusetts, 1994.


\bibitem{GoodeyWeil:1992} P. Goodey and W. Weil,
{\it Centrally symmetric convex bodies and the spherical Radon
 transform}, Journal of Differential Geometry, {\bf 35} (1992), 675--688.


%$\bibitem{Helgason1999} S. Helgason, The Radon Transform, 2nd
%edition, Birkhuser, Boston, 1999.

\bibitem{Helgason:1999}
S. Helgason, {\it The Radon transform}, Volume 5 of progress in
mathematics, Birkhuser Boston Inc., Boston, MA, 1999.


\bibitem{Helgason:2011}
S. Helgason, {\it Integral Geometry and Radon Transforms},
Springer, 2011.


\bibitem{HielscherQuellmalz:2016}
R. Hielscher and M. Quellmalz, {\it  Reconstructing a function on
the sphere from its means along vertical slices}, Inverse Probl.
Imaging, {\bf 10}:3 (2016), 711---739.

\bibitem{HristovaaMoonSteinhauer:2016}
Y. Hristovaa, S. Moon  and D. Steinhauer, {\it A Radon-type
transform arising in Photoacoustic Tomography with circular
detectors: spherical geometry}, Inverse Problems in Science and
Engineering, {\bf 24}:6 (2016), 974--989.
%https://arxiv.org/pdf/1501.02532.pdf.

%S.\,G. Kazantsev  , Th. Schuster, Asymptotic inversion formulas in
%3D vector field tomography for different geometries
%//Journal of Inverse and Ill-posed Problems,    2011, v 19, N  4--5, pp. 769--799.

\bibitem{Kazantsev:2015}
 S.\,G.  Kazantsev, {\it Singular value decomposition for the cone-beam
 transform in the ball},
 J. Inverse Ill-Posed Probl.,
 {\bf  23}:2  (2015),  173--185. %DOI: 10.1515/jiip-2013-0067

\bibitem{LouisRiplingerSpiessSpodarev:2011}
 A.\,K. Louis, M. Riplinger, M. Spiess and E. Spodarev,
 {\it Inversion algorithms for the spherical Radon and cosine transform},
 Inverse Problems, {\bf 27} (2011),  035015. %25 pp.

\bibitem{Louis:2016}
A.\,K. Louis, {\it Exact cone beam reconstruction formulae for
functions and their gradients for spherical and flat detectors},
Inverse Problems, {\bf 32}:11 (2016), 115005. %(16 pp).

\bibitem{Michel:2013}
V. Michel, {\it  Lectures on Constructive Approximation: Fourier,
Spline, and Wavelet Methods on the Real Line, the Sphere, and the
Ball},  Birkhauser, New York, 2013.

\bibitem{Minkowski:1904}
%H. Minkowski, {\it Uber die Korper Konstanter Breite}, Moskau
%Mathematische Sammlung 25 (1904-06), 505--508.
 H. Minkowski, {\it
 About bodies of constant width}, Mathematics Sbornik,  {\bf 25}
  (1904),  505--508. %in Russian.

\bibitem{MorseFeshbach:1953}
F. Morse and  H. Feshbach, {\it  Methods of theoretical physics},
Vol. 2,  McGraw-Hill, 1953.


\bibitem{Muller:1980}
C. Muller, {\it Spherical Harmonics}, Lecture Notes in
Mathematics, Vol. 17, Springer, Berlin, 1966.

\bibitem{NattererWubelling:2001}
F. Natterer and F. Wubelling, {\it Mathematical methods in image reconstruction},
Monographs on Mathematical Modeling and Computation 5,
SIAM, Philadelphia, PA 2001.


\bibitem{Nedelec:2001}
J.-C. Nedelec, {\it Acoustic and Electromagnetic Equations:
Integral Representations for Harmonic Problems},
 Series: Applied Mathematical Sciences, Vol. 144., Springer, New York,
 2001.

%%\bibitem{GelfandShilov:1964}
%%I.\,M. Gelfand and G.\,E. Shilov, {\it Generalized Functions 1,
%%Properties and Operations}, Academic Press, New York-London, 1964.


%[9]  E. Spodarev.  On the rose of intersections for stationary
%flat processes.Adv. Appl. Probab., 33(3):584--599,2001.

%%(3) Inversion algorithms for the spherical Radon and cosine
%transform. Available from:
%https://www.researchgate.net/publication/228931482_Inversion_algorithms_for_the_spherical_Radon_and_cosine_transform
%[accessed Apr 21 2018].

%\bibitem{Helgason1999}
%S. Helgason, The Radon Transform, 2nd edition, Birkhauser, 1999.
%Helgason S. The Radon transform. Progress in mathematics. Boston:
%Birkhauser; 1999.

%\bibitem{Helgason2011}
%S. Helgason, Integral geometry and Radon transforms. New York
%(NY): Springer; 2011.


\bibitem{Palamodov:2017}
V. Palamodov, {\it Reconstruction from cone integral transforms},
  Inverse Problems, {\bf 33}:10 (2017),    104001.




\bibitem{Palamodov:2016}
V.\,P. Palamodov, {\it  Reconstruction from Integral Data},
Monographs and Research Notes in Mathematics. CRC Press, Boca
Raton, 2016.




\bibitem{QuellmalzHielscherLouis:2018}
M. Quellmalz, R. Hielscher, A.\,K. Louis,
 {\it The cone-beam transform and spherical convolution operators},
 (2018), https://arxiv.org/pdf/1803.10515.pdf.

\bibitem{Quellmalz:2017}
M. Quellmalz, {\it  A generalization of the Funk--Radon
transform},  Inverse Problems, {\bf 33}:3 (2017), 035016.

%%!!
\bibitem{Rubin1:1998}
B. Rubin, {\it
Fractional calculus and wavelet transforms in integral geometry},
Fract. Calc. Appl. Anal.,
 {\bf 1}:2 (1998), 193--219.%%!!

%%!!
\bibitem{Rubin2:1998}
B. Rubin, {\it
Spherical Radon transforms and related wavelet transforms},
Appl. Comput. Harmon. Anal., {\bf 5} (1998), 202--215.

%%!!
\bibitem{Rubin:1998}
B. Rubin, {\it Inversion of fractional integrals related to the
spherical Radon transform}, Journal of Functional Analysis, {\bf 157}:2
(1998), 470--487.%%!!


%%!!
\bibitem{Rubin:2000}
 B. Rubin, {\it  Generalized Minkowski--Funk transforms and small
denominators on the sphere}, Fract. Calc. Appl. Anal.,
 {\bf 3}:2 (2000), 177--204.

%%!!
\bibitem{Rubin:2002}
B. Rubin, {\it  Inversion formulas for the spherical Radon
transform and the generalized cosine transform}, Adv. in Appl.
Math., {\bf 29} (2002),  471--497.


%[77] B. Rubin, Inversion formulas for the spherical Radon
%transform and the generalized cosine transform, Advances in
%Applied Mathematics, 29 (2002), pp. 471--97

%%!!
\bibitem{Rubin:2003}
B. Rubin, Notes on Radon transforms in integral geometry, Fract.
Calc. Appl. Anal.,  {\bf 6} (2003), 25--72.


%!!
\bibitem{Rubin:2008}
B. Rubin, {\it  Intersection bodies and generalized cosine
transforms}, Advances in Mathematics, {\bf 218} (2008), 696--727.
%In: htp:/arxiv.org/PS a he/arxiv/pdf/0704/0704.061v2.pdf

%%!!
\bibitem{Rubin:2014}
 B. Rubin, {\it The $\lambda $-cosine transforms with odd kernel
 and the hemispherical transform}, Fract. Calc. Appl. Anal., {\bf 17}:3 (2014),
 765--806.

\bibitem{Rubin:2015}
B. Rubin, {\it Introduction to Radon Transforms: With Elements of
Fractional Calculus and Harmonic Analysis}, Encyclopedia of
Mathematics and its Applications, New York: Cambridge University
Press, 2015.

%\bibitem{Rubin:2016}
%B. Rubin, {\it Radon transforms and Gegenbauer--Chebyshev integrals, }
%Analysis and Mathematical Physics, {\bf 7(2)} (2017) , 117--150.
%\bibitem{Rubin:2016}
% B. Rubin,
% {\it Radon transforms and
% Gegenbauer--Chebyshev integrals, II; examples}, Anal. Math. Phys.,
% Advance online publication, 2016.



\bibitem{Salman:2016}
Y. Salman, {\it  An inversion formula for the spherical transform
in ${\mathbb S}^2$ for a special family of circles of
integration},  Anal. Math. Phys., {\bf 6}:1 (2016), 43--58.


\bibitem{Samko:1983}
S.\,G. Samko, {\it Generalized Riesz potentials and hypersingular
integrals with homogeneous characteristics; their symbols and
inversion}, Trudy Mat. Inst. Steklov.,  {\bf 156} (1980),
157--222.


\bibitem{Samko:2002}
  S.\,G. Samko,
  {\it Hypersingular integrals and their applications},
  Taylor $\&$ Francis, Series: Analytic Methods and Special Functions,
  Vol. 5, 2002.

%\bibitem{SamkoVakulov:2000}
%S.\,G. Samko, B.\,G. Vakulov, {\it On equivalent norms in
%fractional order function spaces of continuous functions on the
%unit sphere}, Fract. Calc. Appl. Anal., {\bf 3:4} (2000),
%401--433.

%\bibitem{Samko1983}
%S. G. Samko,
%Generalized Riesz potentials and hypersingular integrals with
%homogeneous characteristics; their symbols and inversion
%(Russian). Trudy Mat. Inst. Steklov, 156(1980), 157--222 (Transl.
%in Proc. Steklov Inst. Mat. (1983), Issue 2, 173--243) (Math.
%Reviews, 83a: 45004)



%[46] S.\,G. Samko, {\it The Fourier transform of the functions
% $Y_m(x/|x|)/|x|^{n+\alpha}$, Izv. Vyssh. Uchebn. Zaved. Mat. 7 (194) (1978)
%73--78 (in Russian) .

%[47] S.\,G. Samko, {\it eneralized Riesz potentials and hypersingular
%integrals with homogeneous characteristics, their symbols and
%inversion}, Proc. Steklov Inst. Math. 2 (1983) 173--243.
%  [48] S.\,G.
%Samko, [\it Singular integrals over a sphere and the construction of
 %the characteristic from the symbol},
 %Soviet Math. 27 (4) (1983)
 %35--52.



\bibitem{Schneider:1969}
R. Schneider, {\it Functions on a sphere with vanishing integrals
over certain subspheres}, J. Math. Anal. Appl.,
 {\bf 26} (1969), 381--384.



\bibitem{Semyanistyi:1961}
V.\,I. Semyanistyi, {\it Homogeneous functions and some problems
of integral geometry in the spaces of constant curvature}, Dokl.
Akad. Nauk SSSR, {\bf 136}:2 (1961), 288--291.


%V. Semyanisty, {\it Homogeneous functions and some problems of
%integral  geometry in spaces of constant curvature}, Sov. Math.
%Dokl., 2,   (1961), 59--62.

\bibitem{Semyanistyi:1963}
V.\,I. Semyanistyi, {\it Some integral transformations and
integral geometry in an elliptic space}, Tr. Semin. Vektorn.
Tenzorn. Anal., {\bf 12} (1963), 397--441 (in Russian).



\bibitem{Stepanov:2017}
V.\,N. Stepanov, {\it The method of spherical harmonics for
integral transforms on a sphere}, Mathematical Structures and
Modeling, {\bf 2}:42 (2017), 36--48.



\bibitem{Strichartz:1981}
R.\,S. Strichartz, {\it $L_p$ estimates for Radon transforms in
Euclidean and non-Euclidean spaces}, Duke Math. J., {\bf 48}:4
(1981), 699--727.

\bibitem{VarMosKher:1988}
D.\,A. Varshalovich, A.\,N. Moskalev and V.\,K. Khersonskii, {\it
Quantum Theory of Angular Momentum. Irreducible Tensors, Spherical
Harmonics, Vector Coupling Coefficients, $3nj$ Symbols}, World
Scientific Publishing, Teaneck, 1988.

%D. Varshalovich, A. Moskalev, and V. Khersonskii. Quantum Theory
%of Angular Momentum. World Scientific Publishing, Singapore, 1988.



\bibitem{WangLi:2006}
K. Wang, L. Li, {\it Harmonic analysis and approximation on the
unit sphere}, Science Press, Beijing, 2006.


\bibitem{YarmanYazici:2011}
 C.\,E. Yarman and B. Yazici,
 {\it Inversion of the circular averages transform using the Funk
 transform}, Inverse Problems, {\bf 27}:6 (2011), 065001.
 %Volume 27, Number 6, (2011).

\bibitem{ZangerlScherzer:2010}
G. Zangerl and O. Scherzer, {\it Exact reconstruction in
photoacoustic tomography with circular integrating detectors II:
Spherical geometry}, Math. Methods Appl. Sci., {\bf 33}:15 (2010), 1771--1782.


%\bibitem{Strichartz:1983}
%[36] R.\,S. Strichartz, {\it  Analysis of the Laplacian on the
%complete Riemannian manifold}, Journal of Functional Analysis,
% \bf{52(1):48{79}} (1983), ???-???.

%\bibitem{Taylor1981}
% [38] M.\,E. Taylor, {\it Pseudodifferential Operators},
%  Princeton Mathematical Series, Vol. 34. Princeton University Press,
%Princeton, N.J., 1981.

\end{thebibliography}
\end{document}